\documentclass[titlepage,12pt]{article} 
\usepackage{pstricks,pst-plot} 
\usepackage[title]{appendix}
\usepackage{amssymb,amsthm,amsmath} 
\usepackage[a4paper]{geometry}
\usepackage{datetime2}
\usepackage{enumitem}
\usepackage[utf8]{inputenc}
\usepackage{hyperref}

\date{}

\newcommand{\n}{\mathbb{N}}
\newcommand{\z}{\mathbb{Z}}
\newcommand{\q}{\mathbb{Q}}
\newcommand{\re}{\mathbb{R}}

\newcommand{\ep}{\varepsilon}
\renewcommand{\qed}{{\penalty 10000\mbox{$\quad\Box$}}}

\newcommand{\meas}{\operatorname{meas}}
\newcommand{\dint}{\int\!\!\!\!\int}
\newcommand{\Em}{E_{\langle\mu\rangle}}
\newcommand{\ud}{u_{\delta}}
\newcommand{\udhat}{\widehat{u}_{\delta}}
\newcommand{\ad}{a_{\delta}}
\newcommand{\bd}{b_{\delta}}
\newcommand{\cd}{c_{\delta}}
\newcommand{\dd}{d_{\delta}}
\newcommand{\rd}{r_{\delta}}
\newcommand{\kd}{k_{\delta}}
\newcommand{\sd}{S_{\delta}}
\newcommand{\HH}{\mathcal{H}}
\newcommand{\Hd}{H_{\delta,p}}
\newcommand{\Hdhat}{\widehat{H}_{\delta,p}}
\newcommand{\sn}{\mathbb{S}^{d-1}}
\newcommand{\sorth}{\langle\sigma\rangle^{\perp}}
\newcommand{\sz}[1]{#1_{\sigma,z}}

\newcommand{\idp}{\Lambda_{\delta,p}}
\newcommand{\ip}{\Lambda_{0,p}}

\newcommand{\dzp}{\delta\to 0^{+}}
\newcommand{\glim}{\Gamma\mbox{--}\lim}

\newtheorem{thm}{Theorem}[section]

\newtheorem{rmk}[thm]{Remark}
\newtheorem{prop}[thm]{Proposition}
\newtheorem{defn}[thm]{Definition}
\newtheorem{cor}[thm]{Corollary}

\newtheorem{lemma}[thm]{Lemma}

\title{Optimal constants for a non-local approximation of Sobolev norms and total variation}

\author{Authors? Of what?}

\author{Clara Antonucci\vspace{1ex}\\ 
{\normalsize Scuola Normale Superiore} \\
{\normalsize PISA (Italy)}\\
{\normalsize e-mail: \texttt{clara.antonucci@sns.it}}
\and
Massimo Gobbino\vspace{1ex}\\ 
{\normalsize Universit\`a degli Studi di Pisa} \\
{\normalsize PISA (Italy)}\\  
{\normalsize e-mail: \texttt{massimo.gobbino@unipi.it}}
\and
Matteo Migliorini\vspace{1ex}\\ 
{\normalsize Scuola Normale Superiore} \\
{\normalsize PISA (Italy)}\\
{\normalsize e-mail: \texttt{matteo.migliorini@sns.it}}
\and
Nicola Picenni\vspace{1ex}\\ 
{\normalsize Scuola Normale Superiore} \\
{\normalsize PISA (Italy)}\\
{\normalsize e-mail: \texttt{nicola.picenni@sns.it}}
}

\begin{document}
\maketitle
\begin{abstract}

We consider the family of non-local and non-convex functionals proposed and investigated by J.~Bourgain, H.~Brezis and H.-M.~Nguyen in a series of papers of the last decade. It was known that this family of functionals Gamma-converges to a suitable multiple of the Sobolev norm or the total variation, depending on the summability exponent, but the exact constants and the structure of recovery families were still unknown, even in dimension one. 

We prove a Gamma-convergence result with explicit values of the constants in any space dimension. We also show the existence of recovery families consisting of smooth functions with compact support.

The key point is reducing the problem first to dimension one, and then to a finite combinatorial rearrangement inequality.
	
\vspace{6ex}

\noindent{\bf Mathematics Subject Classification 2010 (MSC2010):}
26B30, 46E35.

% 26B30   	Absolutely continuous functions, functions of bounded variation
% 46E35   	Sobolev spaces and other spaces of "smooth'' functions, embedding theorems, trace theorems
\vspace{6ex}

\noindent{\bf Key words:} Gamma-convergence, Sobolev spaces, bounded variation functions, monotone rearrangement, non-local functional, non-convex functional. 

\vspace{6ex}

% \listoftodos

\end{abstract}

%%%%%%%%%%%%%%%%%%%%%
%                   %
%   Inizio lavoro   %
%                   %
%%%%%%%%%%%%%%%%%%%%%
 
\section{Introduction}

Let $p\geq 1$ and $\delta>0$ be real numbers, let $d$ be a positive integer, and let $\Omega\subseteq\re^{d}$ be an open set. For every measurable function $u:\Omega\to\re$ we set
\begin{equation}
\idp(u,\Omega):=\dint_{I(\delta,u,\Omega)}\frac{\delta^{p}}{|y-x|^{d+p}}\,dx\,dy,
\label{defn:idp}
\end{equation}
where
\begin{equation}
I(\delta,u,\Omega):=\left\{(x,y)\in\Omega^{2}:|u(y)-u(x)|>\delta\right\}.
\nonumber
\end{equation}

This family of non-convex and non-local functionals was introduced, motivated and investigated in a series of papers by H.-M.~Nguyen~\cite{2006-JFA-Nguyen,2007-CRAS-Nguyen,2008-JEMS-Nguyen,2011-Duke-Nguyen,2014-JFP-Nguyen}, J.~Bourgain and H.-M.~Nguyen~\cite{2006-CRAS-BouNgu}, H.~Brezis and H.-M.~Nguyen~\cite{2018-AnPDE-BreNgu} (see also~\cite{2015-Lincei-Brezis} and~\cite{2017-CRAS-BreNgu}).

We point out that the dependence on $u$ is just on the integration set. The fixed integrand is divergent on the diagonal $y=x$, and the integration set is closer to the diagonal where the gradient of $u$ is large. This suggests that $\idp(u,\Omega)$ is proportional, in the limit as $\delta\to 0^{+}$, to some norm of the gradient of $u$, and more precisely to the functional
\begin{equation}
\ip(u,\Omega):=\left\{
\begin{array}{l@{\qquad}l}
    \displaystyle\int_{\Omega}|\nabla u(x)|^{p}\,dx  & \mbox{if $p>1$ and $u\in W^{1,p}(\Omega)$,}   \\[2ex]
    \mbox{total variation of $u$ in $\Omega$}  & \mbox{if $p=1$ and $u\in BV(\Omega)$,}  \\[1ex]
    +\infty & \mbox{otherwise.}
\end{array}\right.
\label{defn:ip}
\end{equation}

It is natural to compare the family (\ref{defn:idp}) with the classical approximations of Sobolev or $BV$ norms, based on non-local convex functionals such as
\begin{equation}
G_{\ep,p}(u,\Omega):=\dint_{\Omega}\frac{|u(y)-u(x)|^{p}}{|y-x|^{p}}\rho_{\ep}(|y-x|)\,dx\,dy,
\label{defn:convex}
\end{equation}
where gradients are replaced by finite differences weighted by a suitable family $\rho_{\ep}$ of mollifiers. The idea of approximating integrals of the gradient with double integrals of difference quotients, where all pairs of distinct points interact, has been considered independently by many authors in different contexts. For example, E.~De Giorgi proposed an approximation of this kind to the Mumford-Shah functional in any space dimension, in order to overcome the anisotropy of the discrete approximation~\cite{1995-SIAM-chambolle}. The resulting theory was put into paper in~\cite{ms}, and then extended in~\cite{tesi-mora} to more general free discontinuity problems, and in particular to Sobolev and $BV$ spaces. In the same years, the case of Sobolev and $BV$ norms was considered in details in~\cite{2001-BouBreMir} (see also~\cite{2004-CalcVar-Ponce}). 

The result, as expected, is that the family $G_{\ep,p}(u,\re^{d})$ converges as $\ep\to 0^{+}$ to a suitable multiple of $\ip(u,\re^{d})$, both in the sense of pointwise convergence, and in the sense of De Giorgi's Gamma-convergence. This provides a characterization of Sobolev functions (if $p>1$), and of bounded variation functions (if $p=1$), as those functions for which the pointwise limit or the Gamma-limit is finite.

From the heuristic point of view, the non-convex approximating family (\ref{defn:idp}) seems to follow a different paradigm. Indeed, it has been observed by J.M.~Morel (as quoted at page~4 of the transparencies of the conference~\cite{Brezis:youtube}) that this definition involves some sort of ``vertical slicing'' that evokes the definition of integral \emph{à la Lebesgue}, in contrast to the definition \emph{à la Riemann} that seems closer to the ``horizontal slicing'' of the finite differences in (\ref{defn:convex}).

From the mathematical point of view, the asymptotic behavior of $(\ref{defn:idp})$ exhibits some unexpected features. In order to state the precise results, let us introduce some notation. Let $\sn:=\{\sigma\in\re^{d}:|\sigma|=1\}$ denote the unit sphere in $\re^{d}$. For every $p\geq 1$ we consider the geometric constant
\begin{equation}
G_{d,p}:=\int_{\sn}|\langle v,\sigma\rangle|^{p}\,d\sigma,
\label{defn:knp}
\end{equation}
where $v$ is any element of $\sn$ (of course the value of $G_{d,p}$ does not depend on the choice of $v$), and the integration is intended with respect to the $(d-1)$-dimensional Hausdorff measure. 

The main convergence results obtained so far can be summed up as follows.
\begin{itemize} 
\item \emph{Pointwise convergence for $p>1$.} For every $p>1$ it turns out that
\begin{equation}
\lim_{\dzp}\idp(u,\re^d)=\frac{1}{p}G_{d,p}\,\ip(u,\re^{d})
\qquad
\forall u\in L^{p}(\re^{d}).
\label{thbibl:pointwise}
\end{equation}

\item \emph{Pointwise convergence for $p=1$.} In the case $p=1$ equality (\ref{thbibl:pointwise}) holds true for every $u\in C^{1}_{c}(\re^{d})$, but there do exist functions $u\in W^{1,1}(\re^{d})$ for which the left-hand side is infinite (while of course the right-hand side is finite). A precise characterization of equality cases is still unknown.

\item \emph{Gamma-convergence for every $p\geq 1$.} For every $p\geq 1$ there exists a constant $C_{d,p}$ such that
\begin{equation}
\glim_{\dzp}\idp(u,\re^{d})=\frac{1}{p}G_{d,p}C_{d,p}\,\ip(u,\re^{d})
\qquad
\forall u\in L^{p}(\re^{d}),
\nonumber
\end{equation}
where the Gamma-limit is intended with respect to the usual metric of $L^{p}(\re^{d})$ (but the result would be the same with respect to the convergence in $L^{1}(\re^{d})$ or in measure). Moreover, it was proved that $C_{d,p}\in(0,1)$, namely the Gamma-limit is always nontrivial but different from the pointwise limit. 

\end{itemize}

As a consequence, again one can characterize the Sobolev space $W^{1,p}(\re^{d})$ as the set of functions in $L^{p}(\re^{d})$ for which the pointwise limit or the Gamma-limit are finite. As for $BV(\re^{d})$, in this setting it can be characterized only through the Gamma-limit.

Some problems remained open, and were stated explicitly in~\cite{2011-Duke-Nguyen,2018-AnPDE-BreNgu}.

\begin{itemize}
  \item \emph{Question 1}. What is the exact value of $C_{d,p}$, at least in the case $d=1$?
  \item \emph{Question 2}. Does $C_{d,p}$ depend on $d$?
  \item \emph{Question 3}. Do there exist recovery families made up of continuous functions, or even of functions of class $C^{\infty}$?
\end{itemize}

In this paper we answer these three questions. Concerning question~1 and~2, we prove that $C_{d,p}$ does not depend on $d$, and coincides with the value $C_{p}$ conjectured in~\cite{2007-CRAS-Nguyen} for the one-dimensional case, namely
\begin{equation}
C_{p}:=\left\{
\begin{array}{l@{\qquad}l}
    \displaystyle\frac{1}{p-1}\left(1-\frac{1}{2^{p-1}}\right)  & \mbox{if $p>1$,}   \\[3ex]
    \log 2  & \mbox{if $p=1$.}  
\end{array}\right.
\label{defn:cp}
\end{equation}

Concerning the third question, we prove that smooth recovery families do exist. Our main result is the following.

\begin{thm}[Gamma-convergence]\label{thm:main}

Let us consider the functionals $\idp$ and $\ip$ defined in (\ref{defn:idp}) and (\ref{defn:ip}), respectively.

Then for every positive integer $d$ and every real number $p\geq 1$ it turns out that
\begin{equation}
\glim_{\dzp}\idp(u,\re^{d})=\frac{1}{p}G_{d,p}C_{p}\,\ip(u,\re^{d})
\qquad
\forall u\in L^{p}(\re^{d}),
\nonumber
\end{equation}
where $G_{d,p}$ is the geometric constant defined in (\ref{defn:knp}), and $C_{p}$ is the constant defined in (\ref{defn:cp}). In particular, the following two statements hold true.
\begin{enumerate}
\renewcommand{\labelenumi}{(\arabic{enumi})}

\item \emph{(Liminf inequality)} For every family $\{\ud\}_{\delta>0}\subseteq L^{p}(\re^{d})$, with $\ud\to u$ in $L^{p}(\re^{d})$ as $\delta\to 0^{+}$, it turns out that
\begin{equation}
\liminf_{\dzp}\idp(\ud,\re^{d})\geq\frac{1}{p}G_{d,p}C_{p}\,\ip(u,\re^{d}).
\label{th:main-liminf}
\end{equation}

\item \emph{(Limsup inequality)} For every $u\in L^{p}(\re^{d})$ there exists a family $\{\ud\}_{\delta>0}\subseteq L^{p}(\re^{d})$, with $\ud\to u$ in $L^{p}(\re^{d})$ as $\delta\to 0^{+}$, such that
\begin{equation}
\limsup_{\dzp}\idp(\ud,\re^{d})\leq\frac{1}{p}G_{d,p}C_{p}\,\ip(u,\re^{d}).
\nonumber
\end{equation}
We can also assume that the family $\{\ud\}$ consists of functions of class $C^{\infty}$ with compact support.

\end{enumerate}

\end{thm}

The proof of this result requires a different approach to the problem, which we briefly sketch below. In previous literature the constant $C_{d,p}$ was defined through some sort of cell problem as
\begin{equation}
\frac{1}{p}G_{d,p}C_{d,p}:=\inf\left\{\liminf_{\dzp}\idp\left(\ud,(0,1)^{d}\right):\ud\to u_{0}\mbox{ in }L^{p}\left((0,1)^{d}\right)\right\},
\nonumber
\end{equation}
where $u_{0}(x)=(x_{1}+\ldots+x_{d})/\sqrt{d}$. Unfortunately, this definition is quite implicit and provides no informations on the structure of the families that approach the optimal value. This lack of structure complicates things, in such a way that  just proving that $C_{d,p}>0$ requires extremely delicate estimates (this is the content of~\cite{2006-CRAS-BouNgu}). On the Gamma-limsup side, since $\idp$ is quite sensitive to jumps, what is difficult is glueing together the recovery families corresponding to different slopes, even in the case of a piecewise affine function in dimension one. This requires a delicate surgery near the junctions (see~\cite{2011-Duke-Nguyen}). Finally, as for question~3, difficulties originate from the lack of convexity or continuity of the functionals (\ref{defn:idp}), which do not seem to behave well under convolution or similar smoothing techniques.

The core of our approach consists in proving that $\idp$ in dimension one behaves well under  \emph{vertical $\delta$-segmentation} and \emph{monotone rearrangement}. We refer to Section~\ref{sec:asympt-cost} for the details, but roughly speaking this means that monotone step functions whose values are consecutive integer multiples of $\delta$ are the most efficient way to fill the gap between any two given levels. The argument is purely one-dimensional, and it is carried out in Proposition~\ref{prop:main}. In turn, the proof relies on a discrete combinatorial rearrangement inequality, which we investigate in Theorem~\ref{thm:dino-discrete} under more general assumptions.

\begin{figure}[htbp]
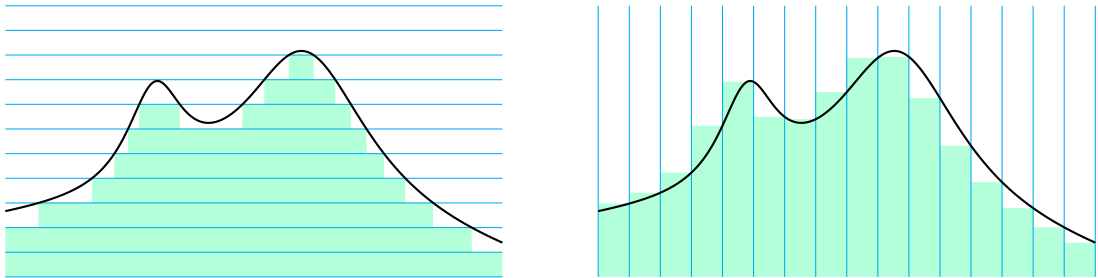

\begin{center}

\SpecialCoor

\newrgbcolor{verdino}{0.7 1 0.85}

\def\funzione#1{6 14 #1 mul sin mul #1 6 sub 2 exp 0.5 mul 1 add div 10 7 #1 mul sin mul #1 12 sub 2 exp 0.1 mul 1 add div add 4 5 #1 mul cos mul add 2 5 #1 mul sin mul add 1.5 div}

\def\Sdelta#1#2{
\pscustom*[linecolor=verdino,plotpoints=2000]{\psplot{0}{20}{#1 #2 div floor #2 mul}\psline(20,0)\psline(0,0)\closepath}
\multido{\n=0+#2}{12}{\psline[linecolor=cyan,linewidth=0.5\pslinewidth](0,\n)(20,\n)}
\psplot[plotpoints=200,linewidth=1\pslinewidth]{0}{20}{#1}
}

\psset{unit=1.8ex}
\hfill
\pspicture(0,-1)(20,11)
\Sdelta{\funzione{x}}{1}
\endpspicture
\hfill\hfill
\pspicture(0,-1)(20,11)
\multido{\n=1.25+1.25}{16}{\psframe*[linecolor=verdino](!\n\space 1.25 sub 0)(!\n\space \funzione{\n})}
\multido{\n=0+1.25}{17}{\psline[linecolor=cyan,linewidth=0.5\pslinewidth](\n,0)(\n,11)}
\psplot[plotpoints=200,linewidth=1\pslinewidth]{0}{20}{\funzione{x}}
\endpspicture
\hfill\mbox{}

\caption{{\small vertical $\delta$-segmentation vs horizontal $\ep$-segmentation ($\delta$ is the distance between the parallel lines on the left, $\ep$ is the distance between the parallel lines on the right)}}
\label{delta-vs-ep}
\end{center}
\end{figure}

We observe that this strategy, namely estimating the asymptotic cost of oscillations by reducing ourselves to a discrete combinatorial minimum problem, is the same exploited in~\cite{ms,tesi-mora}, with the remarkable difference that now the reduction to the discrete setting is achieved through vertical $\delta$-segmentation, while in \cite{ms,tesi-mora} it was obtained through a horizontal $\ep$-segmentation (see Figure~\ref{delta-vs-ep}).

The asymptotic estimate on the cost of oscillations opens the door to the Gamma-liminf inequality in dimension one, which at this point follows from well established techniques. As for the Gamma-limsup inequality, in dimension one we just need to exhibit a family that realizes the given explicit multiple of $\ip(u,\re)$, and this can be achieved through a vertical $\delta$-segmentation \emph{à la Lebesgue} (see Proposition~\ref{prop:limsup}). This produces a recovery family made up of step functions, and it is not difficult to modify them in order to obtain functions of class $C^{\infty}$ with asymptotically the same energy (see Proposition~\ref{prop:d-step}). Finally, passing from dimension one to any dimension is just an application of the one-dimensional result to all the one-dimensional sections of a function of $d$ variables. 

At the end of the day, we have a completely self-contained proof of Theorem~\ref{thm:main} above, and a clear indication that the true difficulty of the problem lies in dimension one, and actually in the discretized combinatorial model. We hope that these ideas could be extended to the more general functionals considered in~\cite{2018-AnPDE-BreNgu}. Some steps in this direction have already been done in~\cite{AGP:log-e} (see also the note~\cite{AGMP:CRAS}).

This paper is organized as follows. In Section~\ref{sec:dino} we develop a theory of monotone rearrangements, first in a discrete, and then in a semi-discrete setting. In Section~\ref{sec:1-dim} we prove our Gamma-convergence result in dimension one. In Section~\ref{sec:n-dim} we prove the Gamma-convergence result in any space dimension.

We would like to thank an anonymous referee for drawing our attention to the rearrangement inequalities of~\cite{1974-AnIF-GarRod}, which could be used in place of our Theorem~\ref{thm:dino-cont} in the appropriate point of the proof of Proposition~\ref{prop:main}. It is not difficult to realize that those rearrangement inequalities, and their discrete combinatorial counterpart known as Taylor's Lemma (see~\cite{1973-JCT-Taylor}), are actually equivalent to ours. Keeping this equivalence into account, we think that our Section~\ref{sec:dino} could be interesting not only because it makes the paper as self-contained as possibile, but also because it provides new proofs of some of the results of~\cite{1973-JCT-Taylor,1974-AnIF-GarRod}, starting from a different perspective (\emph{à la Lebesgue} here, and \emph{à la Riemann} in previous works). It was surprising to discover that similar combinatorial problems originated in the 1970s in a completely different context. 

\setcounter{equation}{0}
\section{An aggregation/segregation problem}\label{sec:dino}

In this section we study the minimum problem for two simplified versions of (\ref{defn:idp}), which we interpret as optimizing the disposition of some objects of different types (actually dinosaurs of different species). The first problem is purely discrete, namely with a finite number of dinosaurs of a finite number of species. The second one is semi-discrete, namely with a continuum of dinosaurs belonging to a finite number of species. 

\subsection{Discrete setting}

Let us consider
\begin{itemize}

\item a positive integer $n$,

\item a function $u:\{1,\ldots,n\}\to\z$,

\item a symmetric subset $E\subseteq\z^{2}$ (namely any subset with the property that $(i,j)\in E$ if and only if $(j,i)\in E$),

\item  a nonincreasing function $h:\{0,1,\ldots,n-1\}\to\re$.
\end{itemize}

Let us introduce the discrete interaction set
\begin{equation}
J(E,u):=\left\{(x,y)\in\{1,\ldots,n\}^{2}: x\leq y,\ (u(x),u(y))\in E\right\},
\label{defn:J-discrete}
\end{equation}
and let us finally define
\begin{equation}
\HH(h,E,u):=\sum_{(x,y)\in J(E,u)}h(y-x).
\label{defn:F-discrete}
\end{equation}

Just to help intuition, we think of $u$ as an arrangement of $n$ dinosaurs placed in the points $\{1,\ldots,n\}$. There are different species of dinosaurs, indexed by integer numbers, so that $u(x)$ denotes the species of the dinosaur in position $x$. The subset $E\subseteq\z^{2}$ is the list of all pairs of species that are hostile to each other. A pair of points $(x,y)$ belongs to $J(E,u)$ if and only if $x\leq y$ and the two dinosaurs placed in $x$ and $y$ belong to hostile species, and in this case the real number $h(y-x)$ measures the ``hostility'' between the two dinosaurs. As expected, the closer are the dinosaurs, the larger is their hostility. 

Keeping this jurassic framework into account, sometimes in the sequel we call $u$ a ``discrete arrangement of $n$ dinosaurs'', we call $E$ an ``enemy list'', we call $h$ a ``discrete hostility function'', and $\HH(h,E,u)$ the ``total hostility of the arrangement''. At this level of generality, we admit the possibility that $(i,i)\in E$ for some integer $i$, namely that a dinosaur is hostile to dinosaurs of the same species, including itself. For this reason, the hostility function $h(x)$ is defined also for $x=0$. This generality turns out to be useful in the proof of the main result for discrete arrangements.

In the sequel we focus on the special case where $E$ coincides with 
\begin{equation}
E_{k}:=\{(i,j)\in\z^{2}:|j-i|\geq k+1\}
\label{defn:E2}
\end{equation}
for some positive integer $k$. In this case it is quite intuitive that the arrangements that minimize the total hostility are the ``monotone'' ones, namely those in which all dinosaurs of the same species are close to each other, and the groups corresponding to different species are sorted in ascending or descending order. To this end, we introduce the following notion.

\begin{defn}[Nondecreasing rearrangement -- Discrete setting]\label{defn:u*-discrete}
\begin{em}

Let $n$ be a positive integer, and let $u:\{1,\ldots,n\}\to\z$ be a function. The \emph{nondecreasing rearrangement} of $u$ is the function $Mu:\{1,\ldots,n\}\to\z$ defined as
\begin{equation}
Mu(x):=\min\left\{j\in\z:|\{y\in\{1,\ldots,n\}:u(y)\leq j\}|\geq x\strut\right\},	
\nonumber
\end{equation}
where $|A|$ denotes the number of elements of the set $A$.
\end{em}
\end{defn}

As the name suggests, $Mu$ turns out to be a nondecreasing function, and it is uniquely characterized by the fact that the two level sets
$$\left\{x\in\{1,\ldots,n\}:u(x)=j\strut\right\},
\qquad\quad
\left\{x\in\{1,\ldots,n\}:Mu(x)=j\strut\right\}$$
have the same number of elements for every $j\in\z$.

As expected, the main result is that monotone arrangements minimize the total hostility with respect to the enemy list $E_{k}$.

\begin{thm}[Total hostility minimization -- Discrete setting]\label{thm:dino-discrete}

Let $n$ and $k$ be two positive integers, let $E_{k}\subseteq\z^{2}$ be the subset defined by (\ref{defn:E2}), and let $h:\{0,\ldots,n-1\}\to\re$ be a nonincreasing function. Let $u:\{1,\ldots,n\}\to\z$ be any function, let $Mu$ be the nondecreasing rearrangement of $u$ introduced in Definition~\ref{defn:u*-discrete}, and let $\HH(h,E_{k},u)$ be the total hostility defined in (\ref{defn:F-discrete}).

Then it turns out that
\begin{equation}
\HH(h,E_{k},u)\geq \HH(h,E_{k},Mu).
\label{th:dino-discrete}
\end{equation}

\end{thm}

\subsection{Semi-discrete setting}

Let us consider
\begin{itemize}

\item  an interval $(a,b)\subseteq\re$,

\item a measurable function $u:(a,b)\to\z$ with finite image,

\item a symmetric subset $E\subseteq\z^{2}$,

\item a nonincreasing function $c:(0,b-a)\to\re$ (note that $c(\sigma)$ might diverge as $\sigma\to 0^{+}$).

\end{itemize}

Let us introduce the semi-discrete interaction set
\begin{equation}
I(E,u):=\left\{(x,y)\in(a,b)^{2}: (u(x),u(y))\in E\right\},
\label{defn:I-continuous}
\end{equation}
and let us finally define
\begin{equation}
\mathcal{F}(c,E,u):=\dint_{I(E,u)}c(|y-x|)\,dx\,dy.
\label{defn:F-continuous}
\end{equation}

In analogy with the discrete setting, we interpret $u(x)$ as a continuous arrangement of dinosaurs of a finite number of species, $c(y-x)$ as the hostility between two dinosaurs of hostile species placed in $x$ and $y$, and we think of $\mathcal{F}(c,E,u)$ as the total hostility of the arrangement $u$ with respect to the enemy list $E$.

Once again, we suspect that monotone arrangements minimize the total hostility with respect to the  enemy list $E_{k}$. This leads to the following notion.

\begin{defn}[Nondecreasing rearrangement -- Semi-discrete setting]\label{defn:u*-continuous}
\begin{em}

Let $u:(a,b)\to\z$ be a measurable function with finite image. The \emph{nondecreasing rearrangement} of $u$ is the function $Mu:(a,b)\to\z$ defined as
\begin{equation}
Mu(x):=\min\left\{j\in\z:\meas\{y\in(a,b):u(y)\leq j\}\geq x-a\strut\right\},	
\nonumber
\end{equation}
where $\meas(A)$ denotes the Lebesgue measure of a subset $A\subseteq(a,b)$.
\end{em}
\end{defn}

The function $Mu$ is nondecreasing and satisfies
$$\meas\{x\in(a,b):u(x)=j\}=\meas\{x\in(a,b):Mu(x)=j\}
\qquad
\forall j\in\z.$$

The following result is the semi-discrete counterpart of Theorem~\ref{thm:dino-discrete}.

\begin{thm}[Total hostility minimization -- Semi-discrete setting]\label{thm:dino-cont}

Let $(a,b)\subseteq\re$ be an interval, let $k$ be a positive integer, let $E_{k}\subseteq\z^{2}$ be the subset defined by (\ref{defn:E2}), and let $c:(0,b-a)\to\re$ be a nonincreasing function. Let $u:(a,b)\to\z$ be any measurable function with finite image, let $Mu$ be the nondecreasing rearrangement of $u$ introduced in Definition~\ref{defn:u*-continuous}, and let $\mathcal{F}(c,E_{k},u)$ be the total hostility defined in (\ref{defn:F-continuous}).

Then it turns out that
\begin{equation}
\mathcal{F}(c,E_{k},u)\geq \mathcal{F}(c,E_{k},Mu).
\label{th:dino-cont}
\end{equation}

\end{thm}

\begin{rmk}\label{rmk:dino-old}
\begin{em}

Theorem~\ref{thm:dino-cont} above is stated in the form that we need in the proof of Proposition~\ref{prop:main}. With a further approximation step in the proof, one can show that the same conclusion (\ref{th:dino-cont}) holds true also without assuming that the image of $u$ is finite and contained in $\z$, and without assuming that $k$ is a positive integer (but just a real number greater than $-1$). 

This extended result is equivalent to~\cite[Theorem~1.1]{1974-AnIF-GarRod}, in the same way as our Theorem~\ref{thm:dino-discrete} is equivalent to Taylor's lemma (see~\cite{1973-JCT-Taylor} and~\cite[Theorem~1.2]{1974-AnIF-GarRod}).

\end{em}
\end{rmk}

\subsection{Proof of Theorem~\ref{thm:dino-discrete}}

Despite the quite intuitive statement, the proof requires some work. To begin with, we introduce some notation, and we develop some preliminary results. Since the function $h$ is fixed once for all, in the sequel we simply write $\HH(E,u)$ instead of $\HH(h,E,u)$.

\paragraph{\textmd{\textit{Left-right gap}}}

Let $L$ and $R$ be two finite sets of positive integers. We call \emph{left-right gap} the quantity
\begin{equation}
\mathcal{G}(L,R):=h(0)+\sum_{\ell\in L}h(\ell)+\sum_{r\in R}h(r)-\sum_{(\ell,r)\in L\times R}\left(h(\ell+r-1)-h(\ell+r)\strut\right).
\label{defn:lr-gap}
\end{equation}

We prove that
\begin{equation}
\mathcal{G}(L,R)\leq\sum_{i=0}^{|L|+|R|}h(i).
\label{th:lr-gap}
\end{equation}

To this end, we argue by induction on the number of elements of $R$. If $R=\emptyset$, from (\ref{defn:lr-gap}) we deduce that
$$\mathcal{G}(L,R):=h(0)+\sum_{\ell\in L}h(\ell)\leq\sum_{i=0}^{|L|}h(i)=\sum_{i=0}^{|L|+|R|}h(i),$$
where the inequality is true term-by-term because $h$ is nonincreasing.

Let us assume now that the conclusion holds true whenever $R$ has $n$ elements, and let us consider any pair $(L,R)$ with $|R|=n+1$. Let us set
$$a:=\max R,
\hspace{4em}
b:=\min\{n\in\n\setminus\{0\}:n\not\in L\},$$
and let us consider the new pair $(L_{1},R_{1})$ defined as
$$L_{1}:=L\cup\{b\},
\qquad\qquad
R_{1}:=R\setminus\{a\}.$$

In words, we have removed the largest element of $R$, and added the smallest possible element to $L$. We observe that $|R_{1}|=n$ and $|L_{1}|+|R_{1}|=|L|+|R|$. Therefore, if we show that
\begin{equation}
\mathcal{G}(L,R)\leq \mathcal{G}(L_{1},R_{1}),
\label{th:l1r1}
\end{equation}
then (\ref{th:lr-gap}) follows from the inductive assumption.

In order to prove (\ref{th:l1r1}), we expand the left-hand side and the right-hand side according to  (\ref{defn:lr-gap}). After canceling out the common terms, with some algebra we obtain that inequality (\ref{th:l1r1}) holds true if and only if
\begin{equation}
h(a)+\sum_{r\in R_{1}}\left(h(b+r-1)-h(b+r)\strut\right)\leq h(b)+\sum_{\ell\in L}\left(h(\ell+a-1)-h(\ell+a)\strut\right).
\label{th:rl}
\end{equation}

All terms in the sums are nonnegative because $h$ is nonincreasing. Let us consider the left-hand side. If $a>1$ we know that $R_{1}\subseteq\{1,\ldots,a-1\}$, and hence
\begin{eqnarray}
h(a)+\sum_{r\in R_{1}}\left(h(b+r-1)-h(b+r)\strut\right) & \leq & h(a)+\sum_{r=1}^{a-1}\left(h(b+r-1)-h(b+r)\strut\right) 
\nonumber  \\
& = & h(a)+h(b)-h(a+b-1).
\label{th:rl-2}
\end{eqnarray}

The same inequality is true for trivial reasons also if $a=1$.

Let us consider now the right-hand side of (\ref{th:rl}). If $b>1$ we know that $L\supseteq\{1,\ldots,b-1\}$, and hence
\begin{eqnarray}
h(b)+\sum_{\ell\in L}\left(h(\ell+a-1)-h(\ell+a)\strut\right) & \geq & h(b)+\sum_{\ell=1}^{b-1}\left(h(\ell+a-1)-h(\ell+a)\strut\right) 
\nonumber \\
& = & h(b)+h(a)-h(a+b-1).
\label{th:rl-1}
\end{eqnarray}

As before, the same inequality is true for trivial reasons also if $b=1$.

Combining (\ref{th:rl-1}) and (\ref{th:rl-2}) we obtain (\ref{th:rl}), which in turn is equivalent to (\ref{th:l1r1}). This completes the proof of (\ref{th:lr-gap}).

\paragraph{\textmd{\textit{Reduction of an arrangement and hostility gap}}}

For every function $v:\{1,\ldots,n\}\to\z$, let us set
$$\mu:=\max\{v(i):i\in\{1,\ldots,n\}\},$$
and let $m\in\{1,\ldots,n\}$ be the \emph{largest index} such that $v(m)=\mu$. 

If $n\geq 2$, we call \emph{reduction} of $v$ the function $Rv:\{1,\ldots,n-1\}\to\z$ defined by
\begin{equation}
[Rv](i):=\left\{\begin{array}{l@{\qquad}l}
     v(i) & \mbox{if }i<m,   \\[0.5ex]
     v(i+1) &   \mbox{if }i\geq m.
\end{array}\right.
\nonumber
\end{equation}

In terms of dinosaurs, $Rv$ is the arrangement obtained from $v$ by removing the rightmost dinosaur of the species indexed by the highest integer, and by shifting one position to the left all subsequent dinosaurs.

When passing from $v$ to $Rv$, the total hostility changes by an amount that we call \emph{hostility gap}, defined as
\begin{equation}
\Delta(E,v):=\HH(E,v)-\HH(E,Rv).
\nonumber
\end{equation}

We observe that interactions between any two dinosaurs placed on the same side of the removed one are equal before and after the removal, and therefore they cancel out when computing the gap. On the contrary, if two hostile dinosaurs are placed within distance $d$ on opposite sides of the removed one, their hostility changes from $h(d)$ to $h(d-1)$ after the removal. It follows that the hostility gap can be written as
\begin{equation}
\Delta(E,v)=\sum_{i\in J_{1}(E,u,m)}h(|m-i|)-\sum_{(i,j)\in J_{2}(E,u,m)}\left(h(j-i-1)-h(j-i)\strut\right),
\label{formula:gap}
\end{equation}
where
\begin{equation}
J_{1}(E,u,m):=\left\{i\in\{1,\ldots,n\}:(u(i),u(m))\in E\strut\right\}
\nonumber
\end{equation}
and
\begin{equation}
J_{2}(E,u,m):=\left\{(i,j)\in\{1,\ldots,n\}^{2}: i<m<j,\ (u(i),u(j))\in E\right\}.
\nonumber
\end{equation}

The first sum in (\ref{formula:gap}) keeps into account the interactions of the removed dinosaur with the rest of the world, the second sum represents the increment of the total hostility due to the reduction of distances among the others.

\paragraph{\textmd{\textit{Monotone rearrangement decreases the hostility gap}}}

We prove that
\begin{equation}
\Delta(E_{k},v)\geq\Delta(E_{k},Mv)
\label{ineq:gap}
\end{equation}
for every arrangement $v$ of $n$ dinosaurs.

In order to prove this inequality, we consider the new enemy list
\begin{equation}
\Em:=\z^{2}\setminus\{\mu,\mu-1,\ldots,\mu-k\}^{2},
\nonumber
\end{equation}
and we claim that
\begin{equation}
\Delta(E_{k},v)\geq\Delta(\Em,v)\geq\Delta(\Em,Mv)=\Delta(E_{k},Mv),
\label{ineq:gap-chain}
\end{equation}
which of course implies (\ref{ineq:gap}).

The equality between the last two terms of (\ref{ineq:gap-chain}) follows from formula (\ref{formula:gap}). Indeed, since $Mv$ is nondecreasing, the removed dinosaur is the rightmost one, and therefore in both cases the second sum in (\ref{formula:gap}) is void. Also the first sum in (\ref{formula:gap}) is the same in both cases, because a dinosaur of the highest species is hostile to another dinosaur with respect to the enemy list $E_{k}$ if and only if it is hostile to the same dinosaur with respect to the enemy list $\Em$.

The inequality between the first two terms of (\ref{ineq:gap-chain}) follows again from formula (\ref{formula:gap}). Indeed, the first sum has the same terms both in the case of the enemy list $E_{k}$ and in the case of the enemy list $\Em$, as observed above. As for the second sum, the interactions with respect to $E_{k}$ are also interactions with respect to $\Em$, and therefore when passing from $E_{k}$ to $\Em$ the second sum cannot decrease. Since the second sum appears in (\ref{formula:gap}) with negative sign, the hostility gap with respect to $\Em$ is less than or equal to the hostility gap with respect to $E_{k}$.

It remains to prove that
\begin{equation}
\Delta(\Em,v)\geq\Delta(\Em,Mv).
\label{ineq:delta-main}
\end{equation}

To this end, we introduce the complement enemy list
$$\Em^{c}:=\{\mu,\mu-1,\ldots,\mu-k\}^{2}=\z^{2}\setminus\Em.$$

Since $\z^{2}$ is the disjoint union of $\Em$ and $\Em^{c}$, and the total hostility is additive with respect to the enemy list, we deduce that
$$\HH(\Em,w)=\HH(\z^{2},w)-\HH(\Em^{c},w)$$
for every arrangement $w$, and for the same reason
$$\Delta(\Em,w)=\Delta(\z^{2},w)-\Delta(\Em^{c},w).$$

Moreover, we observe that the total hostility with respect to $\z^{2}$ depends only on the number of dinosaurs, and in particular $\Delta(\z^{2},v)=\Delta(\z^{2},Mv)$. As a consequence, proving (\ref{ineq:delta-main}) is equivalent to showing that 
\begin{equation}
\Delta(\Em^{c},v)\leq\Delta(\Em^{c},Mv).
\label{ineq:gap-emc}
\end{equation}

The advantage of this ``complement formulation'' is that hostility gaps with respect to $\Em^{c}$ depend only on the relative positions of the removed dinosaur with respect to the other dinosaurs of the species with indices between $\mu-k$ and $\mu$.

To be more precise, let us compute the left-hand side of (\ref{ineq:gap-emc}). Let $m$ denote as usual the position of the dinosaur that is removed from $v$ to $Rv$, and let us set
$$R(v):=\left\{r\geq 1:v(m+r)\in\{\mu,\mu-1,\ldots,\mu-k\}\right\},$$
$$L(v):=\left\{\ell\geq 1:v(m-\ell)\in\{\mu,\mu-1,\ldots,\mu-k\}\right\}.$$

In other words, this means that
$$\{m-\ell:\ell\in L(v)\}\cup\{m\}\cup\{m+r:r\in R(v)\}$$
is the set of all integers $i\in\{1,\ldots,n\}$ such that $v(i)\in\{\mu,\mu-1,\ldots,\mu-k\}$, namely the set of positions where the dinosaurs of the last $k+1$ species are placed.

With this notation, the first sum in (\ref{formula:gap}) is
\begin{equation}
h(0)+\sum_{\ell\in L(v)}h(\ell)+\sum_{r\in R(v)}h(r)
\nonumber
\end{equation}
(we recall that in this ``complement formulation'' the dinosaur in $m$ is also hostile to itself), while the second sum in (\ref{formula:gap}) is
\begin{equation}
\sum_{(\ell,r)\in L(v)\times R(v)}\left(h(\ell+r-1)-h(\ell+r)\strut\right).
\nonumber
\end{equation}

Therefore, it turns out that
\begin{eqnarray}
\Delta(\Em^{c},v) & = & h(0)+\sum_{\ell\in L(v)}h(\ell)+\sum_{r\in R(v)}h(r)-\sum_{(\ell,r)\in L(v)\times R(v)}\left(h(\ell+r-1)-h(\ell+r)\strut\right) 
\nonumber  \\[1ex]
& = & \mathcal{G}(L(v),R(v)),
\label{formula:delta-lhs}
\end{eqnarray}
where $\mathcal{G}$ is the left-right gap defined in (\ref{defn:lr-gap}).

On the other hand, in the nondecreasing arrangement $Mv$ the rightmost dinosaur has $|L(v)|+|R(v)|$ dinosaurs of the last $k+1$ species exactly on its left, and therefore
\begin{equation}
\Delta(\Em^{c},Mv)=\sum_{i=0}^{|L(v)|+|R(v)|}h(i).
\label{formula:delta-rhs}
\end{equation}

Keeping (\ref{formula:delta-lhs}) and (\ref{formula:delta-rhs}) into account, inequality (\ref{ineq:gap-emc}) is exactly (\ref{th:lr-gap}).

\paragraph{\textmd{\textit{Conclusion}}}

We are now ready to prove (\ref{th:dino-discrete}). To this end, we argue by induction on the number $n$ of dinosaurs.

In the case $n=1$ there is nothing to prove.

Let us assume now that the conclusion holds true for all arrangements of $n$ dinosaurs. Let $u$ be any arrangement of $n+1$ dinosaurs, and let $Ru$ denote its reduction. Since $Ru$ is an arrangement of $n$ dinosaurs, from the inductive assumption we know that
\begin{equation}
\HH(E_{k},Ru)\geq \HH(E_{k},M(Ru)).
\nonumber
\end{equation}

On the other hand, from (\ref{ineq:gap}) we know that
\begin{equation}
\Delta(E_{k},u)\geq\Delta(E_{k},Mu).
\nonumber
\end{equation}

Since $M(Ru)=R(Mu)$, we finally conclude that
\begin{eqnarray*}
\HH(E_{k},u) & = & \HH(E_{k},Ru)+\Delta(E_{k},u) \\
 & \geq & \HH(E_{k},M(Ru))+\Delta(E_{k},Mu)  \\
 & = & \HH(E_{k},R(Mu))+\Delta(E_{k},Mu)  \\
 & = & \HH(E_{k},Mu),
\end{eqnarray*}
which proves the conclusion for arrangements of $n+1$ dinosaurs.\qed

\subsection{Proof of Theorem~\ref{thm:dino-cont}}

The proof relies on the following approximation result (we omit the proof, which is an exercise in basic measure theory).

\begin{lemma}\label{lemma:epsilon-sets}

Let $m$ be a positive integer, and let $D_{1}, \ldots,D_{m}$ be disjoint measurable subsets of $(0,1)$ such that
$$\bigcup_{i=1}^{m}D_{i}=(0,1).$$

Then for every $\ep>0$ there exist disjoint subsets $D_{1,\ep}, \ldots, D_{m,\ep}$ of $[0,1]$ such that
$$\bigcup_{i=1}^{m}D_{i,\ep}=(0,1),$$
and such that for every $i=1, \ldots, m$ it turns out that
\begin{itemize}

\item $D_{i,\ep}$ is a finite union of intervals with rational endpoints,

\item the Lebesgue measure of the symmetric difference between $D_{i}$ and $D_{i,\ep}$ is less than or equal to $\ep$.\qed
\end{itemize}

\end{lemma}

We are now ready to prove Theorem~\ref{thm:dino-cont}. First of all, we observe that (\ref{th:dino-cont}) is invariant by translations and homotheties. As a consequence, there is no loss of generality in assuming that $(a,b)=(0,1)$ and $c:(0,1)\to\re$. Then we proceed in three steps. To begin with, we prove (\ref{th:dino-cont}) in the special case where the hostility function $c$ is bounded and the arrangement $u$ has a very rigid structure, then for general $u$ but again bounded hostility function, and finally in the general setting.

\paragraph{\textmd{\textit{Step 1}}}

We prove (\ref{th:dino-cont}) under the additional assumption that the hostility function $c:(0,1)\to\re$ is bounded, and that there exists a positive integer $d$ such that $u(x)$ is constant in each interval of the form $((i-1)/d,i/d)$ with $i=1, \ldots, d$.

Indeed, this is actually the discrete setting. To be more precise, we introduce the discrete arrangement $v:\{1,\ldots,d\}\to\z$ defined as
$$v(i):=u\left(\frac{i-1/2}{d}\right)
\qquad
\forall i\in\{1,\ldots,d\},$$
and the discrete hostility function $h:\{0,\ldots,d-1\}\to\re$ defined as
$$h(i):=\int_{0}^{1/d}dx\int_{i/d}^{(i+1)/d}c(|y-x|)\,dy
\qquad
\forall i\in\{0,\ldots,d-1\},$$
which represents the contribution to the total hostility of two intervals of length $1/d$ occupied by hostile dinosaurs, and placed at distance $i/d$ from each other. Then for every enemy list $E_{k}$ it turns out that
$$\mathcal{F}(c,E_{k},u)=2\HH(h,E_{k},v),$$
where $\HH(h,E_{k},v)$ is the discrete total hostility defined in (\ref{defn:F-discrete}), and the factor~2 keeps into account that both $(x,y)$ and $(y,x)$ are included in the semi-discrete interaction set $I(E_{k},u)$, while only one of them is included in the discrete counterpart $J(E_{k},v)$ (see (\ref{defn:J-discrete}) and (\ref{defn:I-continuous})). Moreover, the monotone rearrangement $Mv$ of $v$ is related to the monotone rearrangement $Mu$ of $u$ by the formula
$$Mv(i)=Mu\left(\frac{i-1/2}{d}\right)
\qquad
\forall i\in\{1,\ldots,d\},$$
and again it turns out that
$$\mathcal{F}(c,E_{k},Mu)=2\HH(h,E_{k},Mv)$$
for every enemy list $E_{k}$. At this point, (\ref{th:dino-cont}) is equivalent to
$$\HH(h,E_{k},v)\geq \HH(h,E_{k},Mv),$$
which in turn is true because of Theorem~\ref{thm:dino-discrete}.

\paragraph{\textmd{\textit{Step 2}}}

We prove (\ref{th:dino-cont}) for a general arrangement $u:(0,1)\to\z$, but again under the additional assumption that the hostility function $c:(0,1)\to\re$ is bounded.

To this end, let $z_{1}<z_{2}<\ldots<z_{m}$ denote the elements in the image of $u$, and let
$$D_{i}:=\{x\in(0,1):u(x)=z_{i}\}
\qquad
\forall i\in\{1,\ldots,m\}$$
denote the set of positions of dinosaurs of the species $z_{i}$. For every $\ep>0$, let us consider the sets $D_{1,\ep}$, \ldots, $D_{m,\ep}$ given by Lemma~\ref{lemma:epsilon-sets}, and the function $u_{\ep}:(0,1)\to\z$ defined as
$$u_{\ep}(x)=z_{i}
\qquad
\forall x\in D_{i,\ep}.$$

Since the hostility function $c$ is bounded, and the symmetric difference between $D_{i}$ and $D_{i,\ep}$ has measure less than or equal to $\ep$, there exists a constant $\Gamma$ (depending on $m$ and $c$, but independent of $\ep$) such that
$$|\mathcal{F}(c,E_{k},u)-\mathcal{F}(c,E_{k},u_{\ep})|\leq \Gamma\ep
\qquad\mbox{and}\qquad
|\mathcal{F}(c,E_{k},Mu)-\mathcal{F}(c,E_{k},Mu_{\ep})|\leq \Gamma\ep.$$

On the other hand, the function $u_{\ep}$ satisfies the assumptions of the previous step, and therefore
$$\mathcal{F}(c,E_{k},u_{\ep})\geq \mathcal{F}(c,E_{k},Mu_{\ep}).$$

From all these inequalities it follows that
$$\mathcal{F}(c,E_{k},u)\geq \mathcal{F}(c,E_{k},Mu)-2\Gamma\ep.$$

Since $\ep>0$ is arbitrary, (\ref{th:dino-cont}) is proved in this case.

\paragraph{\textmd{\textit{Step 3}}}

We prove (\ref{th:dino-cont}) without assuming that the hostility function $c(x)$ is bounded.

To this end, for every $n\in\n$ we consider the truncated hostility function 
$$c_{n}(x):=\min\{c(x),n\}
\qquad
\forall x\in(0,1).$$

We observe that 
$$\mathcal{F}(c,E_{k},u)\geq \mathcal{F}(c_{n},E_{k},u)
\qquad\forall n\in\n$$
because $c(x)\geq c_{n}(x)$ for every $x\in(0,1)$, and
$$\mathcal{F}(c_{n},E_{k},u)\geq \mathcal{F}(c_{n},E_{k},Mu)
\qquad
\forall n\in\n$$
because of the result of the previous step applied to the bounded hostility function $c_{n}(x)$. As a consequence, we obtain that
\begin{equation}
\mathcal{F}(c,E_{k},u)\geq \mathcal{F}(c_{n},E_{k},Mu)
\qquad
\forall n\in\n.
\label{th:dino-cont-n}
\end{equation}

On the other hand, by monotone convergence we deduce that
$$\mathcal{F}(c,E_{k},Mu)=\sup_{n\in\n}\mathcal{F}(c_{n},E_{k},Mu),$$
and therefore (\ref{th:dino-cont}) follows from  (\ref{th:dino-cont-n}).\qed 

\setcounter{equation}{0}
\section{Gamma-convergence in dimension one}\label{sec:1-dim}

In this section we prove Theorem~\ref{thm:main} for $d=1$, in which case \begin{equation}
G_{1,p}=2
\qquad
\forall p\geq 1.
\label{k1p}
\end{equation}

To begin with, we introduce the notion of vertical $\delta$-segmentation, which is going to play a crucial role in many parts of the proof.

\begin{defn}[Vertical $\delta$-segmentation]\label{defn:delta-segm}
\begin{em}

Let $\mathbb{X}$ be any set, let $w:\mathbb{X}\to\re$ be any function, and let $\delta>0$. The vertical $\delta$-segmentation of $w$ is the function $S_{\delta}w:\mathbb{X}\to\re$ defined by
\begin{equation}
S_{\delta}w(x):=\delta\left\lfloor\frac{w(x)}{\delta}\right\rfloor
\qquad
\forall x\in\mathbb{X}.
\label{defn:sd}
\end{equation}

The function $S_{\delta}w$ takes its values in $\delta\z$, and it is uniquely characterized by the fact that $S_{\delta}w(x)=k\delta$ for some $k\in\z$ if and only if $k\delta\leq w(x)<(k+1)\delta$.

\end{em}
\end{defn}

\subsection{Asymptotic cost of oscillations}\label{sec:asympt-cost}

Let us assume that a function $\ud(x)$ oscillates between two values $A$ and $B$ in some interval $(a,b)$. Does this provide an estimate from below for $\idp(\ud,(a,b))$, at least  when $\delta$ is small enough? The following Proposition and the subsequent corollaries give a sharp quantitative answer to this question. They are the fundamental tool in the proof of the liminf inequality.

\begin{prop}[Limit cost of vertical oscillations]\label{prop:main}

Let $p\geq 1$ be a real number, let $(a,b)\subseteq\re$ be an interval, and let $\{\ud\}_{\delta>0}\subseteq L^{p}((a,b))$ be a family of functions.

Let us assume that there exist two real numbers $A\leq B$ such that 
\begin{equation}
\liminf_{\dzp}\meas\{x\in(a,b):\ud(x)\leq A+\ep\}>0
\qquad
\forall\ep>0,
\label{hp:liminf-meas1}
\end{equation}
and
\begin{equation}
\liminf_{\dzp}\meas\{x\in(a,b):\ud(x)\geq B-\ep\}>0
\qquad
\forall\ep>0.
\label{hp:liminf-meas2}
\end{equation}

Then it turns out that
\begin{equation}
\liminf_{\dzp}\idp(\ud,(a,b))\geq\frac{2}{p}\cdot C_{p}\cdot\frac{(B-A)^{p}}{(b-a)^{p-1}},
\label{th:liminf-loc}
\end{equation}
where $C_{p}$ is the constant defined in (\ref{defn:cp}).

\end{prop}

\paragraph{\textmd{\textit{Proof}}}

To begin with, we observe that (\ref{th:liminf-loc}) is trivial if $A=B$, or if the left-hand side is infinite. Up to restricting ourselves to a sequence $\delta_{k}\to 0^{+}$, we can also assume that the liminf is actually a limit. Therefore, in the sequel we assume that the left-hand side of (\ref{th:liminf-loc}) is uniformly bounded from above, and that $A<B$.

Let us fix $\ep>0$ such that $4\ep<B-A$. Due to assumptions~(\ref{hp:liminf-meas1}) and~(\ref{hp:liminf-meas2}), there exist $\eta>0$ and $\delta_{0}>0$ such that
\begin{equation}
\meas\{x\in(a,b):\ud(x)\leq A+\ep\}\geq\eta
\qquad
\forall\delta\in(0,\delta_{0}),
\label{meas-a-1}
\end{equation}
\begin{equation}
\meas\{x\in(a,b):\ud(x)\geq B-\ep\}\geq\eta
\qquad
\forall\delta\in(0,\delta_{0}).
\label{meas-b-1}
\end{equation}

\subparagraph{\textmd{\textit{Truncation, $\delta$-segmentation and monotone rearrangement}}}

In this section of the proof, we replace $\{\ud\}$ with a new family $\{\udhat\}$ of monotone piecewise constant functions that still satisfies (\ref{hp:liminf-meas1}) and (\ref{hp:liminf-meas2}), without increasing the left-hand side of (\ref{th:liminf-loc}). To this end, we perform three operations on $\ud(x)$.

The first operation is a truncation between $A$ and $B$. To be more precise, we define $T_{A,B}\ud:(a,b)\to\re$ by setting
$$T_{A,B}\ud(x):=\left\{
\begin{array}{l@{\qquad}l}
   A   & \mbox{if }\ud(x)<A,   \\[0.5ex]
   \ud(x)   & \mbox{if }A\leq\ud(x)\leq B,   \\[0.5ex]
   B   & \mbox{if }\ud(x)>B.  
\end{array}
\right.$$

We observe that the implication
$$|T_{A,B}\ud(y)-T_{A,B}\ud(x)|>\delta
\quad\Longrightarrow\quad
|\ud(y)-\ud(x)|>\delta$$
holds true for every $x$ and $y$ in $(a,b)$, and hence
$$\idp(T_{A,B}\ud,(a,b))\leq \idp(\ud,(a,b))
\qquad
\forall\delta>0.$$

We also observe that (\ref{meas-a-1}) and (\ref{meas-b-1}) remain true if we replace $\ud(x)$ by $T_{A,B}\ud(x)$.

The second operation is a vertical $\delta$-segmentation, namely we replace $T_{A,B}\ud$ by the function $S_{\delta}T_{A,B}\ud$ defined according to (\ref{defn:sd}).

Again we observe that the implications
\begin{eqnarray*}
|S_{\delta}T_{A,B}\ud(y)-S_{\delta}T_{A,B}\ud(x)|>\delta
 & \Longrightarrow &
|S_{\delta}T_{A,B}\ud(y)-S_{\delta}T_{A,B}\ud(x)|\geq 2\delta \\[1ex]
 & \Longrightarrow & |T_{A,B}\ud(y)-T_{A,B}\ud(x)|>\delta
\end{eqnarray*}
hold true for every $x$ and $y$ in $(a,b)$, and hence
$$\idp(S_{\delta}T_{A,B}\ud,(a,b))\leq \idp(T_{A,B}\ud,(a,b))
\qquad
\forall\delta>0.$$

As for (\ref{meas-a-1}) and (\ref{meas-b-1}), we set $\delta_{1}:=\min\{\ep,\delta_{0}\}$, and we observe that now
\begin{equation}
\meas\{x\in(a,b):S_{\delta}T_{A,B}\ud(x)\leq A+2\ep\}\geq\eta
\qquad
\forall\delta\in(0,\delta_{1}),
\label{meas-a-2}
\end{equation}
\begin{equation}
\meas\{x\in(a,b):S_{\delta}T_{A,B}\ud(x)\geq B-2\ep\}\geq\eta
\qquad
\forall\delta\in(0,\delta_{1}).
\label{meas-b-2}
\end{equation}

The third and last operation we perform is monotone rearrangement, namely we replace $S_{\delta}T_{A,B}\ud$ with the nondecreasing function $MS_{\delta}T_{A,B}\ud$ in $(a,b)$ whose level sets have the same measure of the level sets of $S_{\delta}T_{A,B}\ud$ (see Definition~\ref{defn:u*-continuous}).

From (\ref{meas-a-2}) and (\ref{meas-b-2}) we deduce that now
\begin{equation}
MS_{\delta}T_{A,B}\ud(x)\leq A+2\ep
\qquad
\forall x\in(a,a+\eta),\quad\forall\delta\in(0,\delta_{1}),
\label{meas-a-3}
\end{equation}
\begin{equation}
MS_{\delta}T_{A,B}\ud(x)\geq B-2\ep
\qquad
\forall x\in(b-\eta,b),\quad\forall\delta\in(0,\delta_{1}).
\label{meas-b-3}
\end{equation}

Moreover, we claim that
\begin{equation}
\idp(MS_{\delta}T_{A,B}\ud,(a,b))\leq \idp(S_{\delta}T_{A,B}\ud,(a,b))
\qquad
\forall\delta>0. 
\label{ineq:ud*-sud}
\end{equation}

This is a straightforward consequence of Theorem~\ref{thm:dino-cont}. To be more formal, let us consider the semi-discrete arrangement $v_{\delta}:(a,b)\to\z$ defined by
$$v_{\delta}(x):=\frac{1}{\delta}S_{\delta}T_{A,B}\ud(x)
\qquad
\forall x\in(a,b)$$
(we recall that $S_{\delta}T_{A,B}\ud$ takes its values in $\delta\z$, and hence $v_{\delta}(x)$ is integer valued), and the hostility function $c:(0,b-a)\to\re$ defined as $c(\sigma):=\delta^{p}\sigma^{-1-p}$. We observe that 
$$MS_{\delta}T_{A,B}\ud(x)=\delta Mv_{\delta}(x)
\qquad
\forall x\in(a,b),$$
where $Mv_{\delta}$ is the nondecreasing rearrangement of $v_{\delta}$ according to Definition~\ref{defn:u*-continuous}. 

We observe also that for every pair of points $x$ and $y$ in $(a,b)$ it turns out that
$$(x,y)\in I(\delta,\sd T_{A,B}\ud,(a,b))
\quad\Longleftrightarrow\quad
|v_{\delta}(y)-v_{\delta}(x)|\geq 2
\quad\Longleftrightarrow\quad
(x,y)\in I(E_{1},v_{\delta}),$$
where $E_{1}$ is the enemy list defined in (\ref{defn:E2}), and $I(E_{1},v_{\delta})$ is the semi-discrete interaction set defined according to~(\ref{defn:I-continuous}). It follows that
$$\idp(S_{\delta}T_{A,B}\ud,(a,b))=\mathcal{F}(c,E_{1},v_{\delta}),
\qquad\quad
\idp(MS_{\delta}T_{A,B}\ud,(a,b))=\mathcal{F}(c,E_{1},Mv_{\delta}),$$
and therefore (\ref{ineq:ud*-sud}) is equivalent to (\ref{th:dino-cont}).

In conclusion, the three operations described so far delivered us a family 
$$\udhat:=MS_{\delta}T_{A,B}\ud$$ 
of nondecreasing functions such that the image of $\udhat$ is contained in $\delta\z$. This family satisfies (\ref{meas-a-3}) and (\ref{meas-b-3}), and
\begin{equation}
\idp(\ud,(a,b))\geq \idp(\udhat,(a,b))
\qquad
\forall\delta>0.
\label{ineq:ud*-ud}
\end{equation}

In the sequel we are going to show that any such family satisfies
\begin{equation}
\liminf_{\delta\to 0^{+}}\idp(\udhat,(a,b))\geq\frac{2}{p}\cdot C_{p}\cdot\frac{(B-A-4\ep)^{p}}{(b-a)^{p-1}}.
\label{th:liminf-loc-ep}
\end{equation}
Due to (\ref{ineq:ud*-ud}) and the arbitrariness of $\ep>0$, this is enough to prove (\ref{th:liminf-loc}).

\subparagraph{\textmd{\textit{Extension of the integrals to a vertical strip}}}

In this section of the proof we modify the domain of integration in order to simplify the computation of $\idp(\udhat,(a,b))$.

To begin with, we observe that
$$\idp(\udhat,(a,b))=\dint_{A_{\delta}}\frac{\delta^{p}}{|y-x|^{1+p}}\,dx\,dy\geq\dint_{B_{\delta}}\frac{\delta^{p}}{|y-x|^{1+p}}\,dx\,dy,$$
where
\begin{equation}
A_{\delta}:=I(\delta,\udhat,(a,b))=\left\{(x,y)\in(a,b)^{2}:|\udhat(y)-\udhat(x)|>\delta\right\},
\nonumber
\end{equation}
\begin{equation}
B_{\delta}:=\left\{(x,y)\in(a+\eta,b-\eta)\times(a,b):|\udhat(y)-\udhat(x)|>\delta\strut\right\}.
\nonumber
\end{equation}

Then we write the last integral in the form
$$\dint_{B_{\delta}}\frac{\delta^{p}}{|y-x|^{1+p}}\,dx\,dy=\dint_{B_{\delta}\cup C_{\delta}}\frac{\delta^{p}}{|y-x|^{1+p}}\,dx\,dy-\dint_{C_{\delta}}\frac{\delta^{p}}{|y-x|^{1+p}}\,dx\,dy,$$
where
$$C_{\delta}:=(a+\eta,b-\eta)\times(\re\setminus(a,b)).$$

In other words, the set $B_{\delta}\cup C_{\delta}$ consists of the vertical strip $(a+\eta,b-\eta)\times\re$ minus the set of points $(x,y)\in(a+\eta,b-\eta)\times(a,b)$ such that $|\udhat(y)-\udhat(x)|\leq\delta$. 

Now we observe that
$$\dint_{C_{\delta}}\frac{\delta^{p}}{|y-x|^{1+p}}\,dx\,dy=2\delta^{p}\int_{a+\eta}^{b-\eta}dx\int_{b}^{+\infty}\frac{1}{|y-x|^{1+p}}\,dy.$$

From the convergence of the last double integral it follows that
$$\lim_{\delta\to 0^{+}}\dint_{C_{\delta}}\frac{\delta^{p}}{|y-x|^{1+p}}\,dx\,dy=0,$$
and therefore
\begin{eqnarray}
\liminf_{\delta\to 0^{+}}\idp(\udhat,(a,b))  &  \geq  &  \liminf_{\delta\to 0^{+}}\dint_{B_{\delta}}\frac{\delta^{p}}{|y-x|^{1+p}}\,dx\,dy
\nonumber \\
& = & \liminf_{\delta\to 0^{+}}\dint_{B_{\delta}\cup C_{\delta}}\frac{\delta^{p}}{|y-x|^{1+p}}\,dx\,dy.
\label{est:bd-cd}
\end{eqnarray}

\subparagraph{\textmd{\textit{Computing the integrals}}}

In this last part of the proof we show that
\begin{equation}
\liminf_{\dzp}\dint_{B_{\delta}\cup C_{\delta}}\frac{\delta^{p}}{|y-x|^{1+p}}\,dx\,dy\geq \frac{2}{p}\cdot C_{p}\cdot\frac{(B-A-4\ep)^{p}}{(b-a)^{p-1}}.
\label{th:liminf-final}
\end{equation}

Recalling (\ref{est:bd-cd}), this proves (\ref{th:liminf-loc-ep}), and hence also (\ref{th:liminf-loc}).

To this end, we need to introduce some notation. We know that $\udhat$ is a nondecreasing function with finite image. Let us consider the partition
$$a=x_{0}<x_{1}<\ldots<x_{n}=b$$
of $(a,b)$ with the property that $\udhat(x)$ is constant in each interval of the form $(x_{i-1},x_{i})$, and different intervals correspond to different constants. Let us set
\begin{equation}
h:=\min\{i\in\{1,\ldots,n\}:x_{i}\geq a+\eta\},
\nonumber
\end{equation}
\begin{equation}
k:=\max\{i\in\{0,\ldots,n-1\}:x_{i}\leq b-\eta\}.
\nonumber
\end{equation}

Of course $n$, $h$, $k$, as well as the partition, do depend on $\delta$. Now we claim that
\begin{equation}
\dint_{B_{\delta}\cup C_{\delta}}\frac{\delta^{p}}{|y-x|^{1+p}}\,dx\,dy\geq \frac{2}{p}\cdot C_{p}\cdot\frac{\delta^{p}(k-h-1)^{p}}{(b-a)^{p-1}}
\qquad
\forall\delta\in(0,\delta_{1}).
\label{th:liminf-bdcd}
\end{equation}

To begin with, we show that the values of $\udhat(x)$ in neighboring intervals are consecutive multiples of $\delta$, namely if $\udhat(x)=m\delta$ in $(x_{i-1},x_{i})$ for some $m\in\z$, then $\udhat(x)=(m+1)\delta$ in $(x_{i},x_{i+1})$. Let us assume indeed that $\udhat(x)\geq(m+2)\delta$ in $(x_{i},x_{i+1})$. In this case it turns out that
$$\idp(\udhat,(a,b))\geq\int_{x_{i-1}}^{x_{i}}dx\int_{x_{i}}^{x_{i+1}}\frac{\delta^{p}}{(y-x)^{1+p}}\,dy,$$
which is absurd because the left-hand side is uniformly bounded from above and the integral in the right-hand side is divergent.

With these notations it turns out that
\begin{eqnarray*}
\lefteqn{\hspace{-4em}\dint_{B_{\delta}\cup C_{\delta}}\frac{\delta^{p}}{|y-x|^{1+p}}\,dx\,dy} 
\\[1ex]
\hspace{4em} & \geq & \sum_{i=h+1}^{k-1}\left(\int_{x_{i-1}}^{x_{i}}dx\int_{x_{i+1}}^{+\infty}\frac{\delta^{p}}{|y-x|^{1+p}}\,dy+\int_{x_{i}}^{x_{i+1}}dx\int_{-\infty}^{x_{i-1}}\frac{\delta^{p}}{|y-x|^{1+p}}\,dy\right)  \\[1ex]
& = & \frac{\delta^{p}}{p}\sum_{i=h+1}^{k-1}\left(\int_{x_{i-1}}^{x_{i}}\frac{1}{(x_{i+1}-x)^{p}}\,dx+\int_{x_{i}}^{x_{i+1}}\frac{1}{(x-x_{i-1})^{p}}\,dx\right).
\end{eqnarray*}

Now we distinguish two cases.
\begin{itemize}
\item If $p=1$, computing the integrals we obtain that
$$\dint_{B_{\delta}\cup C_{\delta}}\frac{\delta}{(y-x)^{2}}\,dx\,dy\geq\delta\sum_{i=h+1}^{k-1}\log\left(\frac{x_{i+1}-x_{i-1}}{x_{i+1}-x_{i}}\cdot\frac{x_{i+1}-x_{i-1}}{x_{i}-x_{i-1}}\right).$$

If $\ell_{i}:=x_{i}-x_{i-1}$ denotes the length of the $i$-th interval of the partition, and we apply the inequality between arithmetic and geometric mean, we obtain that
\begin{eqnarray*}
\dint_{B_{\delta}\cup C_{\delta}}\frac{\delta}{(y-x)^{2}}\,dx\,dy  &  \geq  &  \delta\sum_{i=h+1}^{k-1}\log\frac{(\ell_{i}+\ell_{i+1})^{2}}{\ell_{i}\cdot\ell_{i+1}}
\nonumber \\
& \geq & \delta\sum_{i=h+1}^{k-1}\log 4
\nonumber \\
& = & 2\log 2\cdot\delta(k-h-1),
\end{eqnarray*}
which proves (\ref{th:liminf-bdcd}) in this case.

\item  If $p>1$, computing the integrals we obtain that
$$\dint_{B_{\delta}\cup C_{\delta}}\frac{\delta^{p}}{|y-x|^{1+p}}\,dx\,dy\geq\frac{\delta^{p}}{p(p-1)}\sum_{i=h+1}^{k-1}\left(\frac{1}{\ell_{i+1}^{p-1}}+\frac{1}{\ell_{i}^{p-1}}-\frac{2}{(\ell_{i+1}+\ell_{i})^{p-1}}\right),$$
where we set $\ell_{i}:=x_{i}-x_{i-1}$ as before. Therefore, with two applications of Jensen's inequality to the convex function $t\to t^{1-p}$, we obtain that
\begin{eqnarray*}
\dint_{B_{\delta}\cup C_{\delta}}\frac{\delta^{p}}{|y-x|^{1+p}}\,dx\,dy & \geq & \frac{\delta^{p}}{p(p-1)}\sum_{i=h+1}^{k-1}\frac{2^{p}-2}{(\ell_{i+1}+\ell_{i})^{p-1}} 
\nonumber \\[1ex]
 & \geq & \frac{\delta^{p}(2^{p}-2)}{p(p-1)}\cdot\frac{(k-h-1)^{p}}{\left(\sum_{i=h+1}^{k-1}(\ell_{i+1}+\ell_{i})\right)^{p-1}} 
\nonumber \\[1ex]
 & \geq & \frac{\delta^{p}(2^{p}-2)}{p(p-1)}\cdot\frac{(k-h-1)^{p}}{(2(b-a))^{p-1}} 
\nonumber \\[1ex]
 & = & \frac{2}{p}\cdot C_{p}\cdot\frac{\delta^{p}(k-h-1)^{p}}{(b-a)^{p-1}},
\end{eqnarray*}
which proves (\ref{th:liminf-bdcd}) also in this case.
\end{itemize}

Now it remains to estimate $\delta(k-h-1)$. To this end, from (\ref{meas-a-3}) and the minimality of $h$ we deduce that
$$A+2\ep\geq\udhat(x)=:m_{A}\delta
\qquad
\forall x\in(x_{h-1},x_{h}).$$

Similarly, from (\ref{meas-b-3}) and the maximality of $k$ we deduce that
$$B-2\ep\leq\udhat(x)=:m_{B}\delta
\qquad
\forall x\in(x_{k},x_{k+1}).$$

Since the values of $\udhat$ in consecutive intervals are consecutive multiples of $\delta$, it turns out that
$$m_{B}=m_{A}+(k-h+1),$$
and therefore
$$(k-h-1)\delta=(k-h+1)\delta-2\delta=(m_{B}-m_{A})\delta-2\delta\geq B-A-4\ep-2\delta.$$

Plugging this inequality into (\ref{th:liminf-bdcd}), and letting $\delta\to 0^{+}$, we obtain (\ref{th:liminf-final}), which completes the proof.\qed
\bigskip

The following result is a straightforward consequence of Proposition~\ref{prop:main}.

\begin{cor}\label{cor:liminf}

Let us assume that $\ud\to u$ in $L^{p}(\re)$, and let $(a,b)\subseteq\re$ be an interval whose endpoints $a$ and $b$ are Lebesgue points of $u$.

Then it turns out that
\begin{equation}
\liminf_{\dzp}\idp(\ud,(a,b))\geq\frac{2}{p}\cdot C_{p}\cdot\frac{|u(b)-u(a)|^{p}}{(b-a)^{p-1}}.
\nonumber
\end{equation}

\end{cor}

\paragraph{\textmd{\textit{Proof}}}

It is enough to apply Proposition~\ref{prop:main} with $A:=\min\{u(a),u(b)\}$ and $B:=\max\{u(a),u(b)\}$. Assumptions (\ref{hp:liminf-meas1}) and (\ref{hp:liminf-meas2}) are satisfied because $a$ and $b$ are Lebesgue points of the limit of the sequence $\ud$.\qed
\bigskip

We conclude with another variant of Proposition~\ref{prop:main}. We do not need this statement in the sequel, but we think that it clarifies once more the relation between oscillations of $\ud$ and values of $\idp(\ud,(a,b))$. 

\begin{cor}

Let $(a,b)\subseteq\re$ be an interval, let $\{\ud\}_{\delta>0}\subseteq L^{p}((a,b))$ be a family of functions, and let $\operatorname{osc}(\ud,(a,b))$ denote the essential oscillation of $\ud$ in $(a,b)$.

Then it turns out that
\begin{equation}
\liminf_{\dzp}\idp(\ud,(a,b))\geq\frac{2}{p}C_{p}\frac{1}{(b-a)^{p-1}}\left(\liminf_{\dzp}\operatorname{osc}(\ud,(a,b))\right)^{p}.
\nonumber
\end{equation}

\end{cor}

\paragraph{\textmd{\textit{Proof}}}

Let $i_{\delta}$ and $s_{\delta}$ denote the essential infimum and the essential supremum of $\ud(x)$ in $(a,b)$, respectively. Let us assume that $i_{\delta}$ and $s_{\delta}$ are real numbers (otherwise an analogous argument works with standard minor changes). Let us set $w_{\delta}(x):=\ud(x)-i_{\delta}$, and let us observe that
$$\idp(\ud,(a,b))=\idp(w_{\delta},(a,b))
\qquad
\forall \delta>0.$$

Now it is enough to apply Proposition~\ref{prop:main} with $A:=0$ and 
$$B:=\liminf_{\dzp}(s_{\delta}-i_{\delta})=\liminf_{\dzp}\operatorname{osc}(\ud,(a,b)).\qed$$

\subsection{Gamma-liminf}

Corollary~\ref{cor:liminf} represents a ``localized'' version of the liminf inequality (\ref{th:main-liminf}), which now follows from well established techniques (see for example \cite{ms,tesi-mora}). To this end, we need the following characterization of $\ip(u,\re)$ (we omit the standard proof, based on the convexity of the norm).

\begin{lemma}[Piecewise affine horizontal segmentation]\label{lemma:pwa}

Let $p\geq 1$ be a real number, and let $u\in L^{p}(\re)$.

Then there exists $c\in\re$ such that $c+q$ is a Lebesgue point of $u$ for every $q\in\q$. 

Moreover, if for every positive integer $k$ we consider the piecewise affine function $v_{k}:\re\to\re$ such that
\begin{equation}
v_{k}\left(c+\frac{i}{k}\right)=u\left(c+\frac{i}{k}\right)
\qquad
\forall i\in\z,
\nonumber
\end{equation}
then it turns out that
\begin{equation}
\ip(u,\re)=\lim_{k\to +\infty}\int_{\re}|v_{k}'(x)|^{p}\,dx=\sup_{k\geq 1}\int_{\re}|v_{k}'(x)|^{p}\,dx.\qed
\nonumber
\end{equation}

\end{lemma}

We are now ready to prove (\ref{th:main-liminf}) in the case $d=1$. Let $\ud\to u$ be any family converging in $L^{p}(\re)$, and let $c$ and $v_{k}$ be as in Lemma~\ref{lemma:pwa}. For every $i\in\z$, we set $c_{k,i}:=c+i/k$, and we apply Corollary~\ref{cor:liminf} in the interval $(c_{k,i},c_{k,i+1})$. We obtain that
$$\liminf_{\dzp}\idp\left(\ud,(c_{k,i},c_{k,i+1})\right)\geq\frac{2}{p} C_{p}\frac{\left|u(c_{k,i+1})-u(c_{k,i})\right|^{p}}{(1/k)^{p-1}}=\frac{2}{p}C_{p}\int_{c_{k,i}}^{c_{k,i+1}}|v_{k}'(x)|^{p}\,dx.$$

Since
$$\idp(\ud,\re)\geq\sum_{i\in\z}\idp(\ud,(c_{k,i},c_{k,i+1}))
\qquad
\forall\delta>0,$$
we deduce that
\begin{eqnarray*}
\liminf_{\dzp}\idp(\ud,\re) & \geq & \liminf_{\dzp}\sum_{i\in\z}\idp(\ud,(c_{k,i},c_{k,i+1}))\\[1ex]
& \geq & \sum_{i\in\z}\liminf_{\dzp}\idp(\ud,(c_{k,i},c_{k,i+1})) \\[1ex]
& \geq & \frac{2}{p}C_{p}\sum_{i\in\z}\int_{c_{k,i}}^{c_{k,i+1}}|v_{k}'(x)|^{p}\,dx \\[1ex]
& = & \frac{2}{p}C_{p}\int_{\re}|v_{k}'(x)|^{p}\,dx.
\end{eqnarray*}

Letting $k\to +\infty$, and recalling (\ref{k1p}), we obtain exactly (\ref{th:main-liminf}).\qed

\subsection{Gamma-limsup}

This subsection is devoted to a proof of statement~(2) of Theorem~\ref{thm:main} in the case $d=1$.

It is well-known that we can limit ourselves to showing the existence of recovery families for every $u$ belonging to a subset of $L^{p}(\re)$ that is dense in energy with respect to $\ip(u,\re)$. Classical examples of subsets that are dense in energy are the space $C^{\infty}_{c}(\re)$ of functions of class $C^{\infty}$ with compact support, or the space of piecewise affine functions with compact support. Here for the sake of generality we consider the space $PC^{1}_{c}(\re)$ of piecewise $C^{1}$ functions with compact support, defined as follows.

\begin{defn}\label{defn:pc1}
\begin{em}

Let $u:\re\to\re$ be a function. We say that $u\in PC^{1}_{c}(\re)$ if $u$ has compact support, it is Lipschitz continuous, and there exists a \emph{finite} subset $S\subseteq\re$ such that $u\in C^{1}(\re\setminus S)$.

\end{em}
\end{defn}

We show that for every $u\in PC^{1}_{c}(\re)$ the family $\sd u$ of vertical $\delta$-segmentations of $u$ is a recovery family. This proves the Gamma-limsup inequality in dimension one. 

\begin{prop}[Existence of recovery families]\label{prop:limsup}

Let $p\geq 1$ be a real number, and let $u\in PC^{1}_{c}(\re)$ be a piecewise $C^{1}$ function with compact support according to Definition~\ref{defn:pc1}. For every $\delta>0$, let $S_{\delta}u$ denote the vertical $\delta$-segmentation of $u$ according to Definition~\ref{defn:delta-segm}.

Then it turns out that
\begin{equation}
\limsup_{\delta\to 0^{+}}\idp(\sd u,\re)\leq\frac{2}{p}C_{p}\int_{\re}|u'(x)|^{p}\,dx.
\label{th:limsup}
\end{equation}

\end{prop}

\paragraph{\textmd{\textit{Proof}}}

To begin with, we introduce some notation. Let $R_{0}\geq 1$ be any real number such that the support of $u$ is contained in $[-R_{0}+1,R_{0}-1]$. Let $L$ be the Lipschitz constant of $u$ in $\re$, and let $S\subseteq\re$ be a finite set such that $u\in C^{1}(\re\setminus S)$. For every $x\in\re$ and every $\delta>0$ we set
\begin{equation}
J(\delta,u,x):=\{y\in\re:|\sd u(y)-\sd u(x)|>\delta\},
\label{defn:Jdux}
\end{equation}
and
\begin{equation}
\Hd(x):=\int_{J(\delta,u,x)}\frac{\delta^{p}}{|y-x|^{1+p}}\,dy,
\nonumber
\end{equation}
so that
\begin{equation}
\idp(\sd u,\re)=\int_{\re}\Hd(x)\,dx
\qquad
\forall\delta>0.
\label{idp-hdx}
\end{equation}

In the sequel we call $\Hd(x)$ the ``pointwise hostility function''. It represents the contribution of each point $x$ to the double integral defining $\idp(\sd u,\re)$.

\subparagraph{\textmd{\textit{Strategy of the proof}}}

The outline of the proof is the following. First of all, we show that
\begin{equation}
\lim_{\delta\to 0^{+}}\int_{-\infty}^{-R_{0}}\Hd(x)\,dx=\lim_{\delta\to 0^{+}}\int_{R_{0}}^{+\infty}\Hd(x)\,dx=0.
\label{th:limsup-tail}
\end{equation}

Then we define an averaged pointwise hostility function $\Hdhat(x)$ with the property that
\begin{equation}
\int_{-R_{0}}^{R_{0}}\Hd(x)\,dx=\int_{-R_{0}}^{R_{0}}\Hdhat(x)\,dx.
\label{eqn:Hd-Hdhat}
\end{equation}

We also show that the averaged pointwise hostility function satisfies the uniform bound
\begin{equation}
\Hdhat(x)\leq\frac{2}{p}L^{p}
\qquad
\forall x\in[-R_{0},R_{0}],\quad\forall\delta>0,
\label{est:Hdhat}
\end{equation}
and the asymptotic estimate
\begin{equation}
\limsup_{\delta\to 0^{+}}\Hdhat(x)\leq\frac{2}{p}C_{p}|u'(x)|^{p}
\qquad
\forall x\in[-R_{0},R_{0}]\setminus S.
\label{th:limsup-Hdhat}
\end{equation}

At this point, from Fatou's lemma we deduce that
\begin{eqnarray*}
\limsup_{\delta\to 0^{+}}\int_{-R_{0}}^{R_{0}}\Hd(x)\,dx & = & \limsup_{\delta\to 0^{+}}\int_{-R_{0}}^{R_{0}}\Hdhat(x)\,dx \\[1ex]
 & \leq & \int_{-R_{0}}^{R_{0}}\limsup_{\delta\to 0^{+}}\Hdhat(x)\,dx  \\[1ex]
 & \leq & \frac{2}{p}C_{p}\int_{-R_{0}}^{R_{0}}|u'(x)|^{p}\,dx.
\end{eqnarray*}

Keeping (\ref{idp-hdx}) and (\ref{th:limsup-tail}) into account, this estimate implies (\ref{th:limsup}).

\subparagraph{\textmd{\textit{Reducing integration to a bounded interval}}}

We prove (\ref{th:limsup-tail}).

To this end, let us consider any $x\leq -R_{0}$. We observe that in this case the set $J(\delta,u,x)$ defined in (\ref{defn:Jdux}) is contained in the support of $u$, and hence
$$\int_{-\infty}^{-R_{0}}\Hd(x)\,dx\leq\delta^{p}\int_{-\infty}^{-R_{0}}dx\int_{-R_{0}+1}^{R_{0}-1}\frac{1}{|y-x|^{1+p}}\,dy.$$

At this point the first limit in (\ref{th:limsup-tail}) follows from the convergence of the double integral. The proof of the second limit is analogous.

\subparagraph{\textmd{\textit{Uniform bound on the pointwise hostility function}}}

We prove that
\begin{equation}
\Hd(x)\leq\frac{2}{p}L^{p}
\qquad
\forall x\in[-R_{0},R_{0}],\ \forall\delta>0.
\label{est:Hd}
\end{equation}

To this end, we observe that the implication
$$|\sd u(y)-\sd u(x)|>\delta
\quad\Longrightarrow\quad
|u(y)-u(x)|>\delta$$
holds true for every $(x,y)\in\re^{2}$. Since $u$ is Lipschitz continuous, we deduce that
$$|\sd u(y)-\sd u(x)|>\delta
\quad\Longrightarrow\quad
|y-x|\geq\frac{\delta}{L},$$
and hence
$$\Hd(x)\leq\int_{|y-x|\geq\delta/L}\frac{\delta^{p}}{|y-x|^{1+p}}\,dy=2\int_{\delta/L}^{+\infty}\frac{\delta^{p}}{z^{1+p}}\,dz=\frac{2}{p}L^{p},$$
as required.

\subparagraph{\textmd{\textit{Averaged pointwise hostility function}}}

In this part of the proof we introduce the averaged pointwise hostility function. To this end, we consider the open set
$$A(u,\delta):=\{x\in(-R_{0},R_{0}):u(x)\not\in\delta\z\}.$$

A connected component $(a,b)$ of $A(u,\delta)$ is called \emph{monotone} if $[a,b]\cap S=\emptyset$, and $|u'(x)|\geq\delta$ for every $x\in[a,b]$. In this case there exists $k\in\z$ such that $u(a)=k\delta$ and $u(b)=k\delta\pm\delta$, where the sign depends on the sign of $u'(x)$ in $(a,b)$. From the Lipschitz continuity of $u$ we deduce that $A(u,\delta)$ has only a finite number of monotone connected components.

The averaged pointwise hostility function $\Hdhat:\re\to\re$ is defined as
$$\Hdhat(x):=\frac{1}{b-a}\int_{a}^{b}\Hd(s)\,ds$$
if $x\in[a,b)$ for some monotone connected component of $A(\delta,u)$, and $\Hdhat(x):=\Hd(x)$ otherwise.

At this point, inequality (\ref{est:Hdhat}) follows from (\ref{est:Hd}), while (\ref{eqn:Hd-Hdhat}) is true because the integrals of $\Hd(x)$ and $\Hdhat(x)$ are the same both in all monotone connected components, and in the complement set.

\subparagraph{\textmd{\textit{Asymptotic estimate in stationary points}}}

We prove that (\ref{th:limsup-Hdhat}) holds true for every $x\in(-R_{0},R_{0})\setminus S$ with $|u'(x)|=0$.

To begin with, we observe that in this case $x\not\in[a,b)$ for every monotone connected component $(a,b)$ of $A(\delta,u)$ (because $|u'(x)|$ is strictly positive in the closure of every monotone connected component), and therefore $\Hdhat(x)=\Hd(x)$ for every $\delta>0$. 

If $J(\delta,u,x)=\emptyset$ for every $\delta>0$, then $u$ is identically null, and the conclusion is trivial. Otherwise $J(\delta,u,x)\neq\emptyset$ when $\delta$ is small enough. In this case, let $\rd$ be the largest positive real number such that 
$$(x-\rd,x+\rd)\cap J(\delta,u,x)=\emptyset,$$ 
so that 
\begin{equation}
\Hd(x)\leq\int_{-\infty}^{x-\rd}\frac{\delta^{p}}{|y-x|^{1+p}}\,dy+\int_{x+\rd}^{+\infty}\frac{\delta^{p}}{|y-x|^{1+p}}\,dy=\frac{2}{p}\left(\frac{\delta}{\rd}\right)^{p}.
\nonumber
\end{equation}

Let $\delta_{k}\to 0^{+}$ be any sequence such that
\begin{equation}
\limsup_{\delta\to 0^{+}}\frac{\delta}{\rd}=\lim_{k\to +\infty}\frac{\delta_{k}}{r_{\delta_{k}}}.
\label{defn:delta-k}
\end{equation}

Up to subsequences, we can also assume that $r_{\delta_{k}}$ tends to some $r_{0}$. If $r_{0}>0$, then the limit in the right-hand side of (\ref{defn:delta-k}) is~0, which proves (\ref{th:limsup-Hdhat}) in this case. If $r_{0}=0$, then from the maximality of $r_{\delta_{k}}$ we deduce that $|u(x\pm r_{\delta_{k}})-u(x)|=\delta_{k}$ for a suitable choice of the sign, which might depend on $k$. In any case, the limit in the right-hand side of (\ref{defn:delta-k}) turns out to be
$$\lim_{k\to +\infty}\frac{\delta_{k}}{r_{\delta_{k}}}=\lim_{k\to +\infty}\frac{|u(x\pm r_{\delta_{k}})-u(x)|}{r_{\delta_{k}}}=|u'(x)|=0,$$
which proves (\ref{th:limsup-Hdhat}) also in this case.

\subparagraph{\textmd{\textit{Asymptotic estimate in non-stationary points}}}

We prove that (\ref{th:limsup-Hdhat}) holds true for every $x\in(-R_{0},R_{0})\setminus S$ with $|u'(x)|>0$.

Let us assume, without loss of generality, that $u'(x)>0$ (the other case is analogous). Then for every $\delta>0$ small enough it turns out that $x$ lies in the closure of a monotone connected component of $A(\delta,u)$. More precisely, there exist four real numbers $\ad$, $\bd$, $\cd$, $\dd$ with
$$\ad<\bd\leq x<\cd<\dd,$$
and $\kd\in\z$ such that
\begin{equation}
u(\ad)=(\kd-1)\delta,
\quad
u(\bd)=\kd\delta,
\quad
u(\cd)=(\kd+1)\delta,
\quad
u(\dd)=(\kd+2)\delta,
\nonumber
\end{equation}
and
\begin{equation}
u(y)\in((\kd-1)\delta,\kd\delta)
\qquad
\forall y\in(\ad,\bd),
\label{th:u-a-b}
\end{equation}
\begin{equation}
u(y)\in(\kd\delta,(\kd+1)\delta)
\qquad
\forall y\in(\bd,\cd),
\label{th:u-b-c}
\end{equation}
\begin{equation}
u(y)\in((\kd+1)\delta,(\kd+2)\delta)
\qquad
\forall y\in(\cd,\dd).
\label{th:u-c-d}
\end{equation}

We observe that $\ad$, $\bd$, $\cd$, and $\dd$ tend to $x$ as $\delta\to 0^{+}$, and hence
\begin{equation}
\lim_{\delta\to 0^{+}}\frac{\delta}{\bd-\ad}=\lim_{\delta\to 0^{+}}\frac{u(\bd)-u(\ad)}{\bd-\ad}=u'(x).
\label{th:lim-a-b}
\end{equation}

Similarly it turns out that
\begin{equation}
\lim_{\delta\to 0^{+}}\frac{\delta}{\cd-\bd}=\lim_{\delta\to 0^{+}}\frac{\delta}{\dd-\cd}=u'(x),
\label{th:lim-b-c}
\end{equation}
and also
\begin{equation}
\lim_{\delta\to 0^{+}}\frac{\delta}{\cd-\ad}=\lim_{\delta\to 0^{+}}\frac{\delta}{\dd-\bd}=\frac{u'(x)}{2}.
\label{th:lim-a-c}
\end{equation}

From (\ref{th:u-a-b}) through (\ref{th:u-c-d}) we deduce that 
$$J(\delta,u,s)\subseteq(-\infty,\ad]\cup[\dd,+\infty)
\qquad
\forall s\in(\bd,\cd).$$ 

It follows that
$$\Hd(s)\leq\int_{\re\setminus(\ad,\dd)}\frac{\delta^{p}}{|y-s|^{1+p}}\,dy=\frac{\delta^{p}}{p}\left(\frac{1}{(\dd-s)^{p}}+\frac{1}{(s-\ad)^{p}}\right)
\qquad
\forall s\in[\bd,\cd),$$
and hence
\begin{equation}
\Hdhat(x)=\frac{1}{\cd-\bd}\int_{\bd}^{\cd}\Hd(s)\,ds\leq\frac{\delta^{p}}{p}\frac{1}{\cd-\bd}\int_{\bd}^{\cd}\left(\frac{1}{(\dd-s)^{p}}+\frac{1}{(s-\ad)^{p}}\right)
\,ds
\label{est:Hhat}
\end{equation}
for every $x\in[\bd,\cd)$. Now we distinguish two cases.
\begin{itemize}

\item If $p=1$, computing the integrals in (\ref{est:Hhat}) we obtain that
$$\Hdhat(x)\leq\frac{\delta}{\cd-\bd}\log\left(\frac{\dd-\bd}{\delta}\cdot\frac{\delta}{\dd-\cd}\cdot\frac{\cd-\ad}{\delta}\cdot\frac{\delta}{\bd-\ad}\right),$$
and therefore (\ref{th:limsup-Hdhat}) follows from (\ref{th:lim-a-b}) through (\ref{th:lim-a-c}).

\item If $p>1$, computing the integrals in (\ref{est:Hhat}) we obtain that
\begin{eqnarray*}
\Hdhat(x) & \leq & \frac{1}{p(p-1)}\frac{\delta}{\cd-\bd}\cdot \\[1ex]
 & & \mbox{}\cdot\left\{\frac{\delta^{p-1}}{(\dd-\cd)^{p-1}}+\frac{\delta^{p-1}}{(\bd-\ad)^{p-1}}-\frac{\delta^{p-1}}{(\dd-\bd)^{p-1}}-\frac{\delta^{p-1}}{(\cd-\ad)^{p-1}}\right\},
\end{eqnarray*}
and therefore also in this case (\ref{th:limsup-Hdhat}) follows from (\ref{th:lim-a-b}) through (\ref{th:lim-a-c}).

\end{itemize}

This completes the proof.\qed

\subsection{Smooth recovery families}

The aim of this subsection is refining the Gamma-limsup inequality by showing the existence of recovery families consisting of $C^{\infty}$ functions with compact support. To this end, we introduce the following notion.

\begin{defn}[$\delta$-step functions]\label{defn:d-step}
\begin{em}

Let $\delta$ be a positive real number. A function $u:\re\to\re$ is called a \emph{$\delta$-step function} if there exists a positive integer $n$, a $(n+1)$-uple $x_{0}<x_{1}<\ldots<x_{n}$ of real numbers, and $(k_{1}, \ldots, k_{n})\in\z^{n}$ such that
\begin{itemize}
  \item $u(x)=0$ for every $x\in(-\infty,x_{0})\cup(x_{n},+\infty)$,
  \item $u(x)=k_{i}\delta$ in $(x_{i-1},x_{i})$ for every $i=1,\ldots,n$,
  \item $|k_{1}|=|k_{n}|=1$ and $|k_{i}-k_{i-1}|=1$ for every $i=2,\ldots,n$.
\end{itemize}

The values of $u(x)$ for $x\in\{x_{0},x_{1},\ldots,x_{n}\}$ are not relevant (just to fix ideas, we can define $u(x_{i})$ as the maximum between the limit of $u(x)$ as $x\to x_{i}^{+}$ and the limit of $u(x)$ as $x\to x_{i}^{-}$).
\end{em}
\end{defn}

Now we show that, for every \emph{fixed} $\delta>0$, every $\delta$-step function can be approximated in energy by functions of class $C^{\infty}$ with compact support. Roughly speaking, this is possible because the rigid structure of $\delta$-step functions allows to control the effect of convolutions, which otherwise is unpredictable due to the sensitivity of the integration region in (\ref{defn:idp}) to small perturbations.

\begin{prop}[Smooth approximation of $\delta$-step functions]\label{prop:d-step}

Let $\delta>0$ and $p\geq 1$ be real numbers, and let $u:\re\to\re$ be a $\delta$-step function.

Then there exists a family $\{u_{\ep}\}_{\ep>0}\subseteq C^{\infty}_{c}(\re)$ such that
\begin{equation}
\lim_{\ep\to 0^{+}}u_{\ep}= u
\qquad
\mbox{in $L^{p}(\re)$},
\nonumber
\end{equation}
and
\begin{equation}
\lim_{\ep\to 0^{+}}\idp(u_{\ep},\re)=\idp(u,\re).
\nonumber
\end{equation}

\end{prop}

\paragraph{\textmd{\textit{Proof}}}

Let $n$, $x_{i}$ and $k_{i}$ be as in the definition of $\delta$-step functions, and let
\begin{equation}
\tau:=\min\{x_{i}-x_{i-1}:i=1,\ldots,n\}
\nonumber
\end{equation}
be the length of the smallest interval of the partition. We observe that points in neighboring intervals do not contribute to the computation of $\idp(u,\re)$. In particular, if we write as usual
$$\idp(u,\re):=\dint_{I(\delta,u,\re)}\frac{\delta^{p}}{|y-x|^{1+p}}\,dx\,dy,$$
then it turns out that
\begin{equation}
|y-x|\geq\tau
\qquad
\forall(x,y)\in I(\delta,u,\re).
\label{defn:tau}
\end{equation}

Let us fix a mollifier $\rho\in C^{\infty}_{c}(\re)$ with
\begin{itemize}
  \item $\rho(x)\geq 0$ for every $x\in\re$,
  \item $\rho(x)=0$ for every $x\in\re$ with $|x|\geq 1$,
  \item $\int_{\re}\rho(x)\,dx=1$,
\end{itemize}
and let us consider the usual regularization by convolution
\begin{equation}
u_{\ep}(x):=\int_{\re}u(x+\ep y)\rho(y)\,dy.
\nonumber
\end{equation}

It is well-known that $u_{\ep}\in C^{\infty}_{c}(\re)$ for every $\ep>0$, and that for every $p\geq 1$ it turns out that $u_{\ep}\to u$ in $L^{p}(\re)$ as $\ep\to 0^{+}$.

Let us assume that $2\ep<\tau$, let us consider the two open sets
\begin{equation}
A_{\ep}:=\bigcup_{i=0}^{n}(x_{i}-\ep,x_{i}+\ep)\subseteq\re,
\qquad
B_{\ep}:=(A_{\ep}\times\re)\cup(\re\times A_{\ep})\subseteq\re^{2},
\nonumber
\end{equation}
and let us write
\begin{equation}
\idp(u_{\ep},\re)=\dint_{I(\delta,u_{\ep},\re)\cap B_{\ep}}\frac{\delta^{p}}{|y-x|^{1+p}}\,dx\,dy+\dint_{I(\delta,u_{\ep},\re)\setminus B_{\ep}}\frac{\delta^{p}}{|y-x|^{1+p}}\,dx\,dy.
\nonumber
\end{equation}

Since the support of $\rho$ is contained in $[-1,1]$, it turns out that $u_{\ep}(x)=u(x)$ for every $x\in\re\setminus A_{\ep}$. It follows that
$$I(\delta,u_{\ep},\re)\setminus B_{\ep}=I(\delta,u,\re)\setminus B_{\ep},$$
and therefore
$$\lim_{\ep\to 0^{+}}\dint_{I(\delta,u_{\ep},\re)\setminus B_{\ep}}\frac{\delta^{p}}{|y-x|^{1+p}}\,dx\,dy=\lim_{\ep\to 0^{+}}\dint_{I(\delta,u,\re)\setminus B_{\ep}}\frac{\delta^{p}}{|y-x|^{1+p}}\,dx\,dy=\idp(u,\re),$$
where the last equality follows from Lebesgue's dominated convergence theorem because $B_{\ep}$ shrinks to a set of null measure. So it remains to show that
\begin{equation}
\lim_{\ep\to 0^{+}}\dint_{I(\delta,u_{\ep},\re)\cap B_{\ep}}\frac{\delta^{p}}{|y-x|^{1+p}}\,dx\,dy=0.
\label{th:smooth-1}
\end{equation}

To this end, from (\ref{defn:tau}) and the properties of the support of the mollifier, we deduce that now
$$|y-x|\geq\tau-2\ep 
\qquad
\forall (x,y)\in I(\delta,u_{\ep},\re),$$
and therefore
\begin{eqnarray*}
\dint_{I(\delta,u_{\ep},\re)\cap B_{\ep}}\frac{\delta^{p}}{|y-x|^{1+p}}\,dx\,dy & \leq & 2\sum_{i=0}^{n}\int_{x_{i}-\ep}^{x_{i}+\ep}dx\int_{|z|\geq\tau-2\ep}\frac{\delta^{p}}{|z|^{1+p}}\,dz \\
 & \leq & 2\sum_{i=0}^{n}\int_{x_{i}-\ep}^{x_{i}+\ep}\frac{2}{p}\frac{\delta^{p}}{|\tau-2\ep|^{p}}dx \\
 & = & \frac{4}{p}\frac{\delta^{p}}{|\tau-2\ep|^{p}}\cdot 2\ep(n+1),
\end{eqnarray*}
which implies (\ref{th:smooth-1}).\qed
\bigskip

We are now ready to show the existence of smooth recovery families. As usual, it is enough to show the existence of such a family for every $u$ in a subset of $L^{p}(\re)$ which is dense in energy for $\ip(u,\re)$. In this case we consider the space $\mbox{\textit{PA}}_{c}(\re)$ of piecewise affine functions with compact support.

Since piecewise affine functions are piecewise $C^{1}$, we know from Proposition~\ref{prop:limsup} that the family $\sd u$ of vertical $\delta$-segmentations of $u$ is a (non-smooth) recovery family for $u$. The key point is that the vertical $\delta$-segmentation of a piecewise affine function with compact support is a $\delta$-step function according to Definition~\ref{defn:d-step}. Thus from Proposition~\ref{prop:d-step} we deduce the existence of a function $\ud\in C^{\infty}_{c}(\re)$ such that
$$\|\ud-\sd u\|_{L^{p}(\re)}\leq\delta
\qquad\mbox{and}\qquad
\idp(\ud,\re)\leq\idp(\sd u,\re)+\delta$$
for every $\delta>0$. This implies that $\{\ud\}$ is a smooth recovery family for $u$.\qed

\setcounter{equation}{0}
\section{Gamma-convergence in any dimension}\label{sec:n-dim}

It remains to prove Theorem~\ref{thm:main} in any space dimension. This follows from well established sectioning techniques. For every $\sigma\in\sn$, let $\sorth$ denote the hyperplane orthogonal to $\sigma$, namely
$$\sorth:=\{z\in\re^{d}:\langle z,\sigma\rangle=0\}.$$

Given any $u:\re^{d}\to\re$, for every $\sigma\in\sn$ and every $z\in\sorth$, we consider the one-dimensional section $\sz{u}:\re\to\re$ defined as
\begin{equation}
\sz{u}(x):=u(z+\sigma x)
\qquad
\forall x\in\re.
\nonumber
\end{equation}

The main idea is that Sobolev norms, total variation, and functionals such as $\idp$ computed in $u$ are a sort of average of the same quantities computed on the one-dimensional sections $\sz{u}$. The result is the following.

\begin{prop}[Integral-geometric representation]\label{prop:int-geo}

 Let $u:\re^{d}\to\re$ be any measurable function. Let $\idp$ and $\ip$ be the functionals defined in (\ref{defn:idp}) and (\ref{defn:ip}), respectively.

\begin{enumerate}
\renewcommand{\labelenumi}{(\arabic{enumi})}

\item For every $p\geq 1$ it turns out that
\begin{equation}
\int_{\sn}d\sigma\int_{\sorth}\ip(\sz{u},\re)\,dz=G_{d,p}\,\ip(u,\re^{d}),
\nonumber
\end{equation}
where $G_{d,p}$ is the geometric constant defined in (\ref{defn:knp}).

\item  For every $\delta>0$ and every $p\geq 1$ it turns out that
\begin{equation}
\int_{\sn}d\sigma\int_{\sorth}\idp(\sz{u},\re)\,dz=2\idp(u,\re^{d}).\qed
\nonumber
\end{equation}

\end{enumerate}

\end{prop}

We skip the details of the proof of Proposition~\ref{prop:int-geo}, which is a simple application of variable changes in multiple integrals. More generally, for every $\sigma\in\sn$ and every $g\in L^{1}(\re^{d})$ it turns out that
$$\int_{\re^{d}}g(y)\,dy=\int_{\sorth}dz\int_{\re}g(z+\sigma x)\,dx,$$
and this is the main ingredient in the proof of statement~(1).

Similarly, for every $g\in L^{1}(\re^{d}\times\re^{d})$ it turns out that
$$\dint_{\re^{d}\times\re^{d}}g(u,v)\,du\,dv=\frac{1}{2}\int_{\sn}d\sigma\int_{\sorth}dz\dint_{\re\times\re}g(z+\sigma x,z+\sigma y)\cdot|y-x|^{d-1}\,dx\,dy,$$
and this is the main ingredient in the proof of statement~(2).\medskip

We are now ready to prove Theorem~\ref{thm:main}.

\paragraph{\textmd{\textit{Gamma-liminf}}}

Let us assume that $\ud\to u$ in $L^{1}(\re^{d})$. Then for every $\sigma\in\sn$ it turns out that
$$\sz{(\ud)}\to\sz{u}
\qquad
\mbox{in }L^{1}(\re)$$
for almost every $z\in\sorth$. Therefore, from the integral-geometric representations of Proposition~\ref{prop:int-geo}, Fatou's lemma, and the one-dimensional result, we obtain that
\begin{eqnarray*}
\liminf_{\delta\to 0^{+}}\idp(\ud,\re^{d}) & = & \liminf_{\delta\to 0^{+}}\frac{1}{2}\int_{\sn}d\sigma\int_{\sorth}\idp(\sz{(\ud)},\re)\,dz \\[1ex]
 & \geq & \frac{1}{2}\int_{\sn}d\sigma\int_{\sorth}\liminf_{\delta\to 0^{+}}\idp(\sz{(\ud)},\re)\,dz \\[1ex]
 & \geq & \frac{1}{2}\int_{\sn}d\sigma\int_{\sorth}\frac{2}{p}C_{p}\,\ip(\sz{u},\re)\,dz \\[1ex]
 & = & \frac{1}{p}G_{d,p}C_{p}\,\ip(u,\re^{d}).
\end{eqnarray*}

\paragraph{\textmd{\textit{Gamma-limsup}}}

Let $u\in C^{\infty}_{c}(\re^{d})$ be any function with compact support. For every $\delta>0$ we consider the vertical $\delta$-segmentation $S_{\delta}u$ of $u$, and we observe that this operation commutes with the one-dimensional sections, in the sense that
$$\sz{(S_{\delta}u)}=S_{\delta}(\sz{u})
\qquad
\forall\sigma\in\sn,\quad\forall z\in\sorth.$$

Therefore, from the integral-geometric representations of Proposition~\ref{prop:int-geo}, Fatou's lemma, and the one-dimensional result, we obtain that
\begin{eqnarray*}
\limsup_{\delta\to 0^{+}}\idp(S_{\delta}u,\re^{d}) & = & \limsup_{\delta\to 0^{+}}\frac{1}{2}\int_{\sn}d\sigma\int_{\sorth}\idp(\sz{(S_{\delta}u)},\re)\,dz \\[1ex]
 & \leq & \frac{1}{2}\int_{\sn}d\sigma\int_{\sorth}\limsup_{\delta\to 0^{+}}\idp(\sz{(S_{\delta}u)},\re)\,dz \\[1ex]
 & \leq & \frac{1}{2}\int_{\sn}d\sigma\int_{\sorth}\frac{2}{p}C_{p}\,\ip(\sz{u},\re)\,dz \\[1ex]
 & = & \frac{1}{p}G_{d,p}C_{p}\,\ip(u,\re^{d}).
\end{eqnarray*}

The $\delta$-independent bounds on $\idp(\sz{(S_{\delta}u)},\re)$ needed in order to apply Fatou's lemma follow from the Lipschitz continuity of $u$ and the boundedness of its support.

\paragraph{\textmd{\textit{Smooth recovery families}}}

It remains to show the existence of smooth recovery families. The strategy is analogous to the one-dimensional case, and therefore we limit ourselves to outlining the argument, sparing the reader all technicalities.

To begin with, we observe that it is enough to construct smooth recovery families for every $u\in \operatorname{\textit{PA}}_{c}(\re^{d})$. In this case, a non-smooth recovery family is provided by the vertical $\delta$-segmentations $\sd u$ of $u$. On the other hand, vertical $\delta$-segmentations of piecewise affine functions with compact support are $\delta$-step functions, and these functions can be approximated in energy by smooth functions. It follows that for every $\delta>0$ there exists $\ud\in C^{\infty}_{c}(\re^{d})$ such that
$$\|\ud-\sd u\|_{L^{p}(\re^{d})}\leq\delta
\qquad\mbox{and}\qquad
\idp(\ud,\re^{d})\leq\idp(\sd u,\re^{d})+\delta,$$
and therefore $\{\ud\}$ is the required recovery family.

The last approximation step can be proved by convolution as we did in Proposition~\ref{prop:d-step}. To be more precise, a $\delta$-step function in dimension~$d$ is a function $v:\re^{d}\to\re$ with the property that there exist a finite set $P_{1}$, \ldots, $P_{m}$ of disjoint open polytopes (bounded intersections of half-spaces), and integer numbers $k_{1}$, \ldots, $k_{m}$ such that 
\begin{itemize}

  \item $v(x)=k_{i}\delta$ in $P_{i}$ for every $i=1, \ldots,m$,
  
  \item $v(x)=0$ in the open set $P_{0}$ defined as the complement set of the closure of $P_{1}\cup\ldots\cup P_{m}$,
  
  \item $|k_{i}-k_{j}|\leq 1$ whenever the closure of $P_{i}$ intersects the closure of $P_{j}$,
  
  \item $|k_{i}|\leq 1$ whenever the closure of $P_{i}$ intersects the closure of $P_{0}$.
  
\end{itemize}

In words, the level sets of a $\delta$-step function are finite unions of polytopes, and values in adjacent regions differ by $\delta$.

The key point is that for every $\delta$-step function $v$ there exists a positive real number $\tau$ such that
$$(x,y)\in I(\delta,v,\re^{d})
\quad\Longrightarrow\quad
|y-x|\geq\tau.$$

As a consequence, when we define $v_{\ep}$ as the convolution of $v$ with a mollifier whose support is contained in the ball with center in the origin and radius $\ep$, we obtain that
$$(x,y)\in I(\delta,v_{\ep},\re^{d})
\quad\Longrightarrow\quad
|y-x|\geq\tau-2\ep,$$
and at this point the conclusion follows exactly as in the proof of Proposition~\ref{prop:d-step}.\qed

\subsubsection*{\centering Acknowledgments}

The second author has been introduced to this family of non-local functionals by the inspiring talk~\cite{Brezis:youtube} given by H.~Brezis during the congress ``A mathematical tribute to Ennio De Giorgi'', held in Pisa in September 2016 in the 20-th anniversary of his death. The same author had been introduced to non-local approximations of free discontinuity problems by E.~De Giorgi himself in the last year of his life.

We are all deeply grateful to both of them.

%\bibliographystyle{MaxNew}
%\bibliography{../../../BibTeX/Nguyen}

\begin{thebibliography}{10}
\providecommand{\url}[1]{\texttt{#1}}
\providecommand{\urlprefix}{URL }
\providecommand{\selectlanguage}[1]{\relax}
\providecommand{\eprint}[2][]{\url{#2}}

\bibitem{AGMP:CRAS}
\textsc{C.~Antonucci}, \textsc{M.~Gobbino}, \textsc{M.~Migliorini},
  \textsc{N.~Picenni}.
\newblock On the gap between {G}amma-limit and pointwise limit for a non-local
  approximation of the total variation.
\newblock ArXiv:1712.04413.

\bibitem{AGP:log-e}
\textsc{C.~Antonucci}, \textsc{M.~Gobbino}, \textsc{N.~Picenni}.
\newblock On the gap between gamma-limit and pointwise limit for a non-local
  approximation of the total variation.
\newblock ArXiv:1708.01231.

\bibitem{2001-BouBreMir}
\textsc{J.~Bourgain}, \textsc{H.~Brezis}, \textsc{P.~Mironescu}.
\newblock Another look at {S}obolev spaces.
\newblock In \emph{Optimal control and partial differential equations}, pages
  439--455. IOS, Amsterdam, 2001.

\bibitem{2006-CRAS-BouNgu}
\textsc{J.~Bourgain}, \textsc{H.-M. Nguyen}.
\newblock A new characterization of {S}obolev spaces.
\newblock \emph{C. R. Math. Acad. Sci. Paris} \textbf{343} (2006), no.~2,
  75--80.

\bibitem{Brezis:youtube}
\textsc{H.~Brezis}.
\newblock Another triumph for {D}e {G}iorgi’s {G}amma convergence.
\newblock \urlprefix\url{https://www.youtube.com/watch?v=1Y6fvZX1fx8}.
\newblock Conference held during the congress ``A mathematical tribute to Ennio
  De Giorgi'' (Pisa, September 2016).

\bibitem{2015-Lincei-Brezis}
\textsc{H.~Brezis}.
\newblock New approximations of the total variation and filters in imaging.
\newblock \emph{Atti Accad. Naz. Lincei Rend. Lincei Mat. Appl.} \textbf{26}
  (2015), no.~2, 223--240.

\bibitem{2017-CRAS-BreNgu}
\textsc{H.~Brezis}, \textsc{H.-M. Nguyen}.
\newblock Non-convex, non-local functionals converging to the total variation.
\newblock \emph{C. R. Math. Acad. Sci. Paris} \textbf{355} (2017), no.~1,
  24--27.

\bibitem{2018-AnPDE-BreNgu}
\textsc{H.~Brezis}, \textsc{H.-M. Nguyen}.
\newblock Non-local {F}unctionals {R}elated to the {T}otal {V}ariation and
  {C}onnections with {I}mage {P}rocessing.
\newblock \emph{Ann. PDE} \textbf{4} (2018), no.~1, 4:9.

\bibitem{1995-SIAM-chambolle}
\textsc{A.~Chambolle}.
\newblock Image segmentation by variational methods: {M}umford and {S}hah
  functional and the discrete approximations.
\newblock \emph{SIAM J. Appl. Math.} \textbf{55} (1995), no.~3, 827--863.

\bibitem{1974-AnIF-GarRod}
\textsc{A.~M. Garsia}, \textsc{E.~Rodemich}.
\newblock Monotonicity of certain functionals under rearrangement.
\newblock \emph{Ann. Inst. Fourier (Grenoble)} \textbf{24} (1974), no.~2, vi,
  67--116.
%\newblock Colloque International sur les Processus Gaussiens et les
%  Distributions Al\'eatoires (Colloque Internat. du CNRS, No. 222, Strasbourg,
%  1973).

\bibitem{ms}
\textsc{M.~Gobbino}.
\newblock Finite difference approximation of the {M}umford-{S}hah functional.
\newblock \emph{Comm. Pure Appl. Math.} \textbf{51} (1998), no.~2, 197--228.

\bibitem{tesi-mora}
\textsc{M.~Gobbino}, \textsc{M.~G. Mora}.
\newblock Finite-difference approximation of free-discontinuity problems.
\newblock \emph{Proc. Roy. Soc. Edinburgh Sect. A} \textbf{131} (2001), no.~3,
  567--595.

\bibitem{2006-JFA-Nguyen}
\textsc{H.-M. Nguyen}.
\newblock Some new characterizations of {S}obolev spaces.
\newblock \emph{J. Funct. Anal.} \textbf{237} (2006), no.~2, 689--720.

\bibitem{2007-CRAS-Nguyen}
\textsc{H.-M. Nguyen}.
\newblock {$\Gamma$}-convergence and {S}obolev norms.
\newblock \emph{C. R. Math. Acad. Sci. Paris} \textbf{345} (2007), no.~12,
  679--684.

\bibitem{2008-JEMS-Nguyen}
\textsc{H.-M. Nguyen}.
\newblock Further characterizations of {S}obolev spaces.
\newblock \emph{J. Eur. Math. Soc. (JEMS)} \textbf{10} (2008), no.~1, 191--229.

\bibitem{2011-Duke-Nguyen}
\textsc{H.-M. Nguyen}.
\newblock {$\Gamma$}-convergence, {S}obolev norms, and {BV} functions.
\newblock \emph{Duke Math. J.} \textbf{157} (2011), no.~3, 495--533.

\bibitem{2014-JFP-Nguyen}
\textsc{H.-M. Nguyen}.
\newblock Estimates for the topological degree and related topics.
\newblock \emph{J. Fixed Point Theory Appl.} \textbf{15} (2014), no.~1,
  185--215.

\bibitem{2004-CalcVar-Ponce}
\textsc{A.~C. Ponce}.
\newblock A new approach to {S}obolev spaces and connections to
  {$\Gamma$}-convergence.
\newblock \emph{Calc. Var. Partial Differential Equations} \textbf{19} (2004),
  no.~3, 229--255.

\bibitem{1973-JCT-Taylor}
\textsc{H.~Taylor}.
\newblock Rearrangements of incidence tables.
\newblock \emph{J. Combinatorial Theory Ser. A} \textbf{14} (1973), 30--36.

\end{thebibliography}

\label{NumeroPagine}

\end{document}